\documentclass[12pt,reqno]{amsart}
\usepackage{amsthm}
\usepackage{cases}
\usepackage{amscd}
\usepackage{amsfonts}

\usepackage{makecell}
\usepackage{amssymb}
\usepackage{amsmath}
\usepackage[dvipdfmx]{graphicx}
\usepackage{booktabs}
\usepackage{amsthm}

\usepackage{subfigure}
\usepackage{algorithm}
\usepackage{algpseudocode}
\renewcommand{\algorithmicrequire}

\usepackage{cleveref}
\usepackage{color}
\usepackage{datetime}

\allowdisplaybreaks

\makeatletter
\newcommand\figcaption{\def\@captype{figure}\caption}
\newcommand\tabcaption{\def\@captype{table}\caption}
\makeatother

\usepackage{geometry}
\usepackage{algorithm}
\usepackage{multirow}
\usepackage{mathrsfs}
\usepackage{ulem}

\newtheorem{theorem}{Theorem}[section]          
\newtheorem{example}[theorem]{Example}          
\newcommand{\beq}{\begin{equation}}
\newcommand{\eeq}{\end{equation}}
\newcommand{\bea}{\begin{eqnarray}}
\newcommand{\eea}{\end{eqnarray}}
\newcommand{\beas}{\begin{eqnarray*}}
\newcommand{\eeas}{\end{eqnarray*}}

\graphicspath{{Figures/}{logo/}}

\begin{document}
\title[IG-PINNs for elliptic interface problems]{IG-PINNs: Interface-gated physics-informed neural networks for solving elliptic interface problems}

\author[J. Zheng, Y. Huang, N. Yi]{Jiachun Zheng$^\dagger$, Yunqing Huang$^{\ddagger,*}$ and Nianyu Yi$^{\S,*}$}

\address{$^\ddagger$ School of Mathematics and Computational Science, Xiangtan University, Xiangtan 411105, P.R.China }
\email{jczheng2022@126.com}

\address{$\ddagger$ National Center for Applied Mathematics in Hunan, Key Laboratory of Intelligent Computing \& Information Processing of Ministry of Education, Xiangtan University, Xiangtan 411105, Hunan, P.R.China} \email{huangyq@xtu.edu.cn}
\address{$^\S$ Hunan Key Laboratory for Computation and Simulation in Science and Engineering; School of Mathematics and Computational Science, Xiangtan University, Xiangtan 411105, P.R.China}
\email{yinianyu@xtu.edu.cn}

\subjclass{}


\begin{abstract}
 In this work, we develop interface-gated physics-informed neural networks (IG-PINNs) to solve elliptic interface equations. In IG-PINNs, we use a fully connected neural network to capture the smooth behavior across the entire domain. In each subdomain separated by the interface, an interface-gated network is utilized to provide corrections at the interface. In the architectural design of the interface-gated network, we introduce a gating mechanism and a level-set function derived from the interface. This design enables the interface-gated network to effectively handle discontinuous jumps across the interface. Some numerical experiments have confirmed the effectiveness of the IG-PINNs, demonstrating higher accuracy compared with PINNs, interface PINNs (I-PINNs) and multi-domain PINNs (M-PINNs). 
\end{abstract}

\keywords{Physics-informed neural networks; Partial differential equations; Gating mechanism; Elliptic interface problems}

\thanks{$^*$ Corresponding author.}

\maketitle


\section{Introduction}\label{section-1}
The solutions of complex physical problems involve discontinuities caused by multiphysics coupling \cite{1,2}, particularly in interface problems. Interface problems are prevalent in various engineering fields, including materials science \cite{3}, fluid dynamics \cite{4}, heat transfer \cite{5, 6} and multiphase flow in porous media \cite{7}, etc. Classical numerical methods have achieved notable success in solving interface problems, such as finite element method \cite{8,9,10}, finite volume method \cite{11} and immersed boundary method \cite{12}, ghost fluid method \cite{13,14} and its extension \cite{15}, matched interface and boundary method \cite{16}, discontinuous Galerkin method \cite{17}, weak Galerkin method \cite{18}, etc. In \cite{19}, Gibou et al. present a review on numerical methods to simulate multiphase and free surface flows, including the ghost fluid method \cite{14} and the Voronoi interface method \cite{20}. And they also proposed several potential future research priorities regarding the use of deep learning to enhance multiphase flow research and simulation. In addition, there are many numerical algorithms for problems on irregular domains \cite{20,21,22,23,24}.

Based on the universal approximation theorem \cite{25}, deep neural networks \cite{26,27,28,29,30} have introduced novel approaches to develop numerical algorithms to solve PDEs. Unlike classical numerical methods, neural networks provide a nonlinear approximation through diverse network architectures and activation functions in hidden layers. This universal approximation property extends conventional approximation techniques from linear space to nonlinear space. One of the most representative neural network methodologies is physics-informed neural networks (PINNs) \cite{31,32,33}. In PINNs framework, the PDEs along with initial and boundary conditions are incorporated into the construction of loss function. By optimizing the loss function through an appropriately chosen optimization algorithm, the optimal network parameters can be obtained. Furthermore, numerous variants of PINNs have been proposed, including physics-informed attention-based neural networks \cite{34}, self-adaptive physics-informed neural networks \cite{35}, uncertainty quantification physics-informed neural networks with adversarial training \cite{36} and extended physics-informed neural networks (XPINNs) \cite{37} based on domain decomposition.

In recent years, PINNs have been successfully implemented in solving various types of PDEs, including stochastic PDEs \cite{38}, fractional-order differential equations \cite{39}, integro-differential equations \cite{40} and interface problems \cite{41}, etc. In \cite{42}, He et al. utilized a least-squares method to reformulate the interface problem and used a deep neural network as an approximation for the solution. Additionally, Larios-C\'{a}rdenas et al. \cite{43} proposed an efficient hybrid machine learning framework, which utilizes a local, data-driven neural network to dynamically correct interface trajectory errors in coarse-grid level set simulations, significantly reducing numerical dissipation and lowering computational costs.
In \cite{44}, Hu et al. transform interface problems into high-dimensional continuous problems by introducing an additional spatial dimension and then construct a neural network with a single hidden layer to approximate the solution. Furthermore, Tseng et al. \cite{45} proposed a novel approach, namely cusp-capturing PINNs, to solve elliptic interface problems. The key idea is to introduce the interface information as additional input into the network. After optimization, the network learns the intrinsic properties of the solution near the interface. In \cite{46}, interface PINNs (I-PINNs) are developed to handle interface problems. I-PINNs use neural networks with the same weights and structure in each subdomain, differing only in activation functions. This significantly reduces the number of trainable parameters compared with XPINNs \cite{37}. Based on I-PINNs, Roy et al. proposed adaptive interface-PINNs (AdaI-PINNs) for solving interface problems \cite{54}. Compared to I-PINNs, AdaI-PINNs eliminate the need for predefined activation functions.  Wu et al. \cite{55} proposed  hard-constraint and domain-decomposition-based physics-informed neural networks (HCD-PINNs) for nonhomogeneous transient thermal analysis. In addition, Tseng et al. \cite{50} proposed discontinuity and cusp-capturing physics-informed neural network (PINNs) for solving the Stokes equations with piecewise constant viscosity and singular forces along the interface. Li et al. \cite{47} proposed two innovative discontinuity-removing PINNs to solve elliptic interface problems with variable coefficients on curved surfaces. In the Eulerian coordinate framework, surface differential operators are converted into conventional differential operators using level-set functions associated with the surface. And the solution is decomposed into continuous and discontinuous components, each represented by a single neural network that can be trained independently or jointly.

A hybrid neural network approach combined with traditional numerical algorithms, termed localized randomized neural networks based on the finite difference method (LRaNN-FDM) \cite{48}, is proposed to solve interface problems. In LRaNN-FDM, a single randomized neural network approximates the solution in each subdomain. Instead of relying on optimization, LRaNN-DM converts the interface problem into a linear system by using finite difference schemes at randomly sampled points, and then solves the system via least-squares methods. In \cite{58}, Mistani et al. proposed the neural bootstrapping method (NBM). This method trains neural network surrogate models by evaluating the finite discretization residuals of the PDEs obtained on implicit Cartesian grids centered on random collocation points. Its advantage lies in the ability to not only leverage the inherent accuracy and convergence properties of advanced numerical methods but also achieve better scalability to higher-order PDEs by restricting optimization to first-order automatic differentiation. Additionally, several other hybrid neural-classical methods have been successfully applied to interface problems, including physics-informed neural networks-Lattice Boltzmann method (PINN-LBM) \cite{56} and high-order hybrid approach integrating neural networks and fast Poisson solvers for elliptic interface problem \cite{49}.

Traditional numerical algorithms and hybrid methods based on mesh discretization have achieved notable success in solving interface problems. However, these methods often face difficulties in handling interface conditions with complex geometries and dynamic boundaries. Moreover, they struggle with high-dimensional problems. In this work, we propose interface-gated physics-informed neural networks (IG-PINNs) to solve elliptic interface problems. IG-PINNs consist of a fully connected neural network and multiple interface-gated neural networks. Firstly, a fully connected neural network is used to learn the smooth behavior over the entire domain. Secondly, the entire domain is divided into several subdomains by interfaces. In each subdomain, an interface-gated network with a gating mechanism \cite{59} and a level-set function derived from the interface is employed to provide corrections at the interface. Finally, the fully connected neural network is combined with the interface-gated networks in each subdomain to obtain a neural network representation for interface problems. The proposed method adopts a domain decomposition strategy similar to I-PINNs \cite{46} and X-PINNs \cite{37}, but differs in the network models used. Although IG-PINNs have more parameters than I-PINNs and X-PINNs, the former provides a better approximation.

The rest of the paper is organized as follows: In section \ref{section-2}, we explain the elliptic interface problems and the IG-PINNs in detail. In section \ref{section-3}, the different elliptic interface problems are used to test the approximation ability of IG-PINNs. Finally, the advantages and inadequacies of IG-PINNs and its further development are summarized in section \ref{section-4}.

\section{Interface-gated PINNs method for interface problem}\label{section-2}
Let $\Omega=\Omega_{1}\bigcup\Omega_{2}\bigcup\Gamma$, where $\Omega\in R^{2}$, $\Gamma$ denotes the interface and divides the whole domain $\Omega$ into $\Omega_{1}$ and $\Omega_{2}$. $\textbf{n}$ denotes the normal vector pointing from $\Omega_{1}$ to $\Omega_{2}$ along the interface $\Gamma$.  We consider the following elliptic interface problems:
\begin{equation}\label{2.1}
\left\{\begin{aligned}
\nabla\cdot(\alpha(x,y)\nabla u(x,y))=&f(x,y),\qquad (x,y)\in\Omega_{1}\cup\Omega_{2},\\
([\![\alpha\nabla u]\!]\cdot \textbf{n}) (x,y)=&\rho(x,y),\qquad (x,y)\in\Gamma,\\
[\![u]\!](x,y)=&\beta(x,y),\qquad (x,y)\in\Gamma,\\
u(x,y)=&h(x,y),\qquad (x,y)\in\partial\Omega,
\end{aligned}\right.
\end{equation}
where
$[\![v]\!]|_{\Gamma}=v|_{\Omega_1}-v|_{\Omega_2}$ stands for the jump of a function $v$ across the interface $\Gamma$, and the diffusion coefficient
$\alpha(x,y)$ is piecewise defined
\begin{equation}\label{2.2}
\begin{split}
\alpha(x,y)=\left\{
	\begin{aligned}
	\alpha_{1}(x,y), \qquad (x,y)\in \Omega_{1},\\
	\alpha_{2}(x,y), \qquad (x,y)\in \Omega_{2}.\\
	\end{aligned}
	\right.
\end{split}
\end{equation}
Note that $u(x,y)$ is continuous in $\Omega_{1}$ and $ \Omega_{2}$, respectively, but jumps across the interface $\Gamma$:
\begin{equation}\label{2.3}
\begin{split}
u(x,y)=\left\{
	\begin{aligned}
	u_{1}(x,y), \qquad (x,y)\in \Omega_{1},\\
	u_{2}(x,y), \qquad (x,y)\in \Omega_{2}.\\
	\end{aligned}
	\right.
\end{split}
\end{equation}

A single neural network struggles to capture jumps at interfaces due to its inherent continuity. Assigning one neural network to each subdomain separated by the interface can effectively capture these jumps. However, in interface problems, local solutions are coupled at the interface. If the interface conditions are not accurately approximated during training, errors can propagate from the interface into the neighboring subdomains. A natural idea is to integrate interface information into the network structure. This prompts the model to focus more on constraints at the interface.

\begin{figure}[H]
    \centering
    \includegraphics[width=14cm,height=11cm]{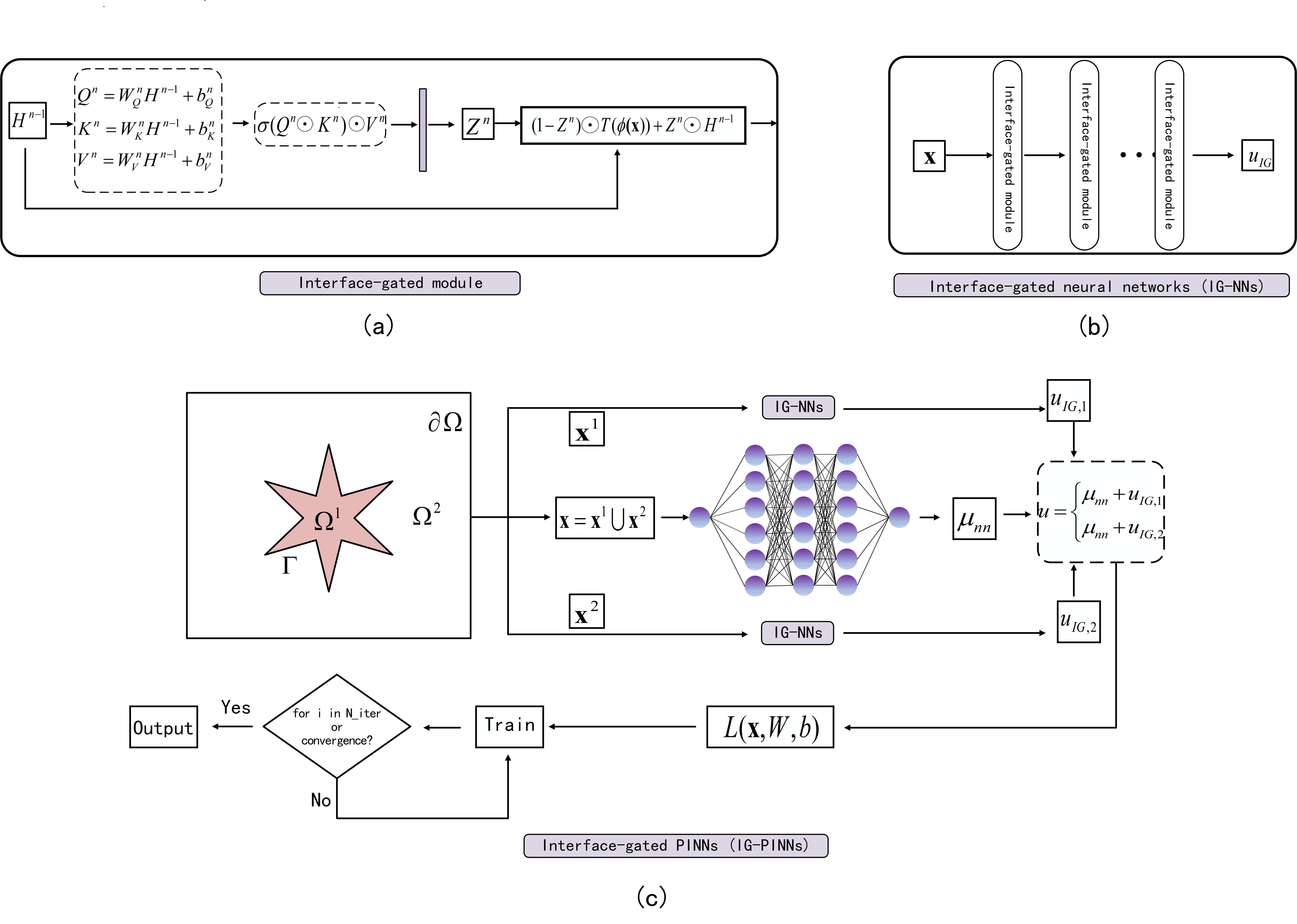}
    \caption{(\textbf{a}): interface-gated module, where the gray rectangle represents a fully connected neural network layer $Z^{n}$. (\textbf{b}): interface-gated neural network (IG-NNs). (\textbf{c}): Interface-gated PINNs (IG-PINNs).}
    \label{F1}
    \end{figure}

In this work, we propose interface-gated physics-informed neural networks (IG-PINNs) for solving elliptic interface problems. IG-PINNs consist of a fully connected neural network and multiple interface-gated neural networks (IG-NNs), where the IG-NNs are composed of several modules with gating mechanisms \cite{59} and level functions derived from interfaces. Within each subdomain separated by the interface, IG-NNs are employed to provide corrections at the interface. It should be noted that within each module of the IG-NNs, the gating mechanism controls the contributions of both internal features and interface-related features. When larger weights are assigned to interface-related features, the model's capacity to learn smooth behaviors within each subdomain separated by the interface may be diminished. Consequently, we need to use a fully connected neural network to learn smooth behaviors across the entire domain. Figure \ref{F1}(a), (b) show the detailed network structure of the IG-NNs. The forward propagation rule for IG-NNs with $M$ module is as follows:
\begin{equation}\label{2.6}
\begin{split}
&T(\phi(\textbf{x}))=\sigma(W^{T}\phi(\textbf{x})+b^{T}),\\
&Z^{n}=\sigma(W^{n}A^{n}+b^{n}),\\
&H^{n}=(1-Z^{n})\odot T(\phi(\textbf{x})) + Z^{n} \odot H^{n-1}, \quad n=1,2,\cdot\cdot\cdot,M,\\
&u_{IG}(\textbf{x},W,b)=W^{M+1}H^{M}+b^{M+1},
\end{split}
\end{equation}
where $H^{0}=\textbf{x}$. $u_{IG}(\textbf{x},W,b)$ is the output of IG-NNs. $\sigma(\cdot)$ denotes the activation function. An adaptive scheme for activation function was proposed in \cite{51} to accelerate convergence in PINNs. $\odot$ denotes a element-wise multiplication. $T(\phi(\textbf{x}))$ is a transmitter network that processes, and transmits the interface information, where $\phi(\textbf{x})$ denotes the level set function derived from the interface, such that $\Gamma=\{\textbf{x} \in R^{d}:\phi(\textbf{x})=0\}$. And $A^{N}$ is defined as follows:
\begin{equation}\label{2.6.1}
\begin{split}
&Q^{n}=W^{n}_{Q}H^{n-1} + b_{Q}^{n},\\
&K^{n}=W^{n}_{K}H^{n-1} + b_{K}^{n},\\
&V^{n}=W^{n}_{V}H^{n-1} + b_{V}^{n}, \\
&A^{n} = \sigma({Q^{n}\odot K^{n}})\odot V^{n}.\\
\end{split}
\end{equation}
IG-NNs consist of multiple modules, where each module contains a gating layer \cite{59} with a residual connection $H^{n-1}$. Here, $Q^{n}, K^{n}$ and $V^{n}$ are used to learn different features and are integrated in the hidden layer $Z^{n}$. $Z^{n}$ and $(1 - Z^{n})$ function as the forget gate and input gate, respectively. Specifically, $Z^{n} \odot H^{n-1}$ controls how much information from the previous module's output is retained, while $(1 - Z^{n}) \odot T(\phi(\textbf{x}))$ controls how much external feature information is forgotten. The expression for a fully connected neural network containing $N$ hidden layers is defined:
\begin{equation}\label{2.7}
\begin{aligned}
&\Phi^{0}(\textbf{x})=\textbf{x},\\
&\Phi^{n}(\textbf{x})=\sigma(W^{n}\Phi^{n-1}(\textbf{x})+b^{n}),\qquad 1\leq n\leq N,\\
&\mu_{nn}(\textbf{x},W,b)=\Phi^{N+1}(\textbf{x})=W^{N+1}\Phi^{N}(\textbf{x})+b^{N+1},\\
\end{aligned}
\end{equation}
where $W^{n}\in R^{m_{n}\times m_{n-1}}$, $b^{n}\in R^{m_{n}}$, $\mu_{nn}(\textbf{x},W,b)$ is the output of the fully connected neural network and $m_{N+1}$ is the output dimension, $\{W, b\}=\{(W^{n},b^{n})_{n=1}^{N+1}\}$ denote all the trainable weights and biases. Thus, a neural network representation for the interface problem (\ref{2.1}) is defined as follows:
\begin{equation}\label{2.8}
u_{nn}(x,y,W, b)=\left\{
	\begin{aligned}
	\mu_{nn}(x,y,W_{\mu}, b_{\mu})+u_{IG,1}(x,y, W_{1}, b_{1}), \qquad (x,y)\in \Omega_{1},\\
	\mu_{nn}(x,y,W_{\mu}, b_{\mu})+u_{IG,2}(x,y, W_{2}, b_{2}), \qquad (x,y)\in \Omega_{2}.\\
	\end{aligned}
	\right.
\end{equation}
 Note that $\mu_{nn}(x,y)$ is eliminated during the calculation of the interface conditions:
\begin{equation}\label{UIG}
\begin{split}
[\![u_{nn}]\!](x_{\Gamma},y_{\Gamma})-\beta(x_{\Gamma},y_{\Gamma})&=(\mu_{nn}(x_{\Gamma},y_{\Gamma})+u_{IG,1}(x_{\Gamma},y_{\Gamma}))-(\mu_{nn}(x_{\Gamma},y_{\Gamma})+u_{IG,2}(x_{\Gamma},y_{\Gamma}))-\beta(x_{\Gamma},y_{\Gamma})\\
&=(u_{IG,1}(x_{\Gamma},y_{\Gamma})-u_{IG,2}(x_{\Gamma},y_{\Gamma}))-\beta(x_{\Gamma},y_{\Gamma}).\\
\end{split}
\end{equation}
Thus, the interface problem (\ref{2.1}) is solved by minimizing the loss function $L(x,y,W,b)$:
\begin{equation}\label{2.9}
\begin{aligned}
L(x,y,W,b)=&\tau\frac{1} {N_{I}}\sum_{i=1}^{N_{I}}| \nabla\cdot(\alpha(x_{i}, y_{i})\nabla u_{nn}(x_{i}, y_{i}, W, b))-f(x_{i}, y_{i})|^{2}\\
&+\tau_{b}\frac{1} {N_{b}}\sum_{i=1}^{N_{b}}|u_{nn}(x_{b}^{i}, y_{b}^{i},W,b)-h(x_{b}^{i}, y_{b}^{i})|^{2}\\
&+\tau_{\Gamma}^{1}\frac{1} {N_{\Gamma}^{1}}\sum_{i=1}^{N_{\Gamma}^{1}}|[\![u_{IG}(x_{\Gamma}^{i}, y_{\Gamma}^{i},W,b)]\!]-\beta(x_{\Gamma}^{i},y_{\Gamma}^{i})|^{2}\\
&+\tau_{\Gamma}^{2}\frac{1} {N_{\Gamma}^{2}}\sum_{i=1}^{N_{\Gamma}^{2}}|[\![\alpha(x_{\Gamma}^{i},y_{\Gamma}^{i}) u_{IG}(x_{\Gamma}^{i},y_{\Gamma}^{i},W,b)]\!]\cdot \textbf{n}-\rho(x_{\Gamma}^{i},y_{\Gamma}^{i})|^{2},\\
\end{aligned}
\end{equation}
where $\{x_{i},y_{i}\}_{i=1}^{N}$, $\{x_{b}^{i},y_{b}^{i}\}_{i=1}^{N_{b}}$ and $\{x_{\Gamma}^{i},y_{\Gamma}^{i}\}_{i=1}^{N_{\Gamma}^{j}}, \;j=1,2$ represent the interior points, the boundary points and the interface points, respectively. The parameters $\{\tau, \tau_{b}, \tau_{\Gamma}^{1}, \tau_{\Gamma}^{2}\}$ represent the weights of each component in the $L(x,y,W,b)$. By minimizing the loss function:
$$\{W^{*},b^{*}\}=\arg \min_{\{W, b\}} L(x,y,W,b),$$
we obtain an approximation $u_{nn}(x,y,W^{*},b^{*})$ of the interface problem (\ref{2.1}). In order to obtain the optimal parameters $\{W^{*},b^{*}\}$, we use gradient descent method to update the parameters. Concretely, the update rules for step $k$ are as follows:
\begin{equation}\label{10}
\begin{aligned}
W^{k+1}=&W^{k}-\eta^{k} \frac{\partial L(x,y,W^{k},b^{k})} {\partial W^{k}},\\
b^{k+1}=&b^{k}-\eta^{k} \frac{\partial L(x,y,W^{k},b^{k})} {\partial b^{k}},\\
\end{aligned}
\end{equation}
where $\eta^{k}$ denotes the step size of the $k$-th iteration. Figure \ref{F1}(c) illustrates the IG-PINNs method in detail.

\section{Numerical experiments }\label{section-3}
In this section, we numerically investigate the performance of the proposed IG-PINNs for solving the elliptic interface equations. For all numerical examples, we set the parameters $\tau=\tau_{b}=\tau_{\Gamma}^{1}=\tau_{\Gamma}^{2}=1$, and chose the stochastic gradient descent method (Adam) \cite{52} with a learning rate 0.001 as the optimizer. The code accompanying this manuscript are publicly available at https://github.com/jczheng126/Interface-gated-PINNs. To better measure the performance of the method, we use the maximum absolute error and the relative $L^{2}$ error:
\begin{equation}\label{3.11}
\begin{split}
&E_{M}= \max_{1\leq i\leq N_{s}}|u(\textbf{x}_{i})-u_{nn}(\textbf{x}_{i})|,\\
&E_{L^2}= \frac{\sqrt{\sum_{i=1}^{N_{s}}|u(\textbf{x}_{i})-u_{nn}(\textbf{x}_{i})|^{2}}} {\sqrt{\sum_{i=1}^{N_{s}}|u(\textbf{x}_{i})|^{2}}} .\\
\end{split}
\end{equation}
Additionally, we use the metrics defined as follows to evaluate the approximation of IG-PINNs and other neural network methods at the interface:

\begin{equation}\label{inter}
\begin{split}
&E_{jump}= \frac{||[\![u_{nn}]\!](\textbf{x})-\beta(\textbf{x})||_{L^{2}(\Gamma)}} {||\beta(\textbf{x})||_{L^{2}(\Gamma)}} ,\\
&E_{flux}= \frac{||([\![\alpha\nabla u_{nn}]\!]\cdot \textbf{n}) (\textbf{x})-\rho(\textbf{x})||_{L^{2}(\Gamma)}} {||\rho(\textbf{x})||_{L^{2}(\Gamma)}} .\\
\end{split}
\end{equation}
To strictly quantify the contribution of each component in IG-PINNs, we conducted comprehensive ablation experiments. These experiments were designed to systematically remove individual components in order to evaluate their importance in isolation. We begin with the ablation experiment by considering the following elliptic interface problem:
\begin{equation}\label{A1}
\left\{\begin{aligned}
\Delta u(x,y)=f(x,y), &\qquad(x,y)\in \Omega/\Gamma,\\
u(x,y)=m(x,y), &\qquad(x,y)\in\partial\Omega,\\
[\![u]\!](x,y)=g(x,y), &\qquad(x,y)\in\Gamma,\\
([\![\nabla u]\!]\cdot\textbf{n})(x,y)=h(x,y),&\qquad(x,y)\in\Gamma,\\
\end{aligned}\right.
\end{equation}
where $\Omega$ is an annular region with an inner radius $r_{inn}=0.151$ and an outer radius $r_{out}=0.913$. The exact solution:
\begin{equation}\label{A2}
u(x,y)=\left\{
	\begin{aligned}
	\sum_{i=1}^{4}\frac{sin(i\pi x)sin(i\pi y)}{4}, \qquad (x,y)\in \Omega_{1},\\
	\sum_{i=1}^{2}\frac{sin(i\pi x)sin(i\pi y)}{2}, \qquad (x,y)\in \Omega_{2}.\\
	\end{aligned}
	\right.
\end{equation}
The interface is described by the following level set function:
$$\phi(x,y)=\sqrt{x^{2}+y^{2}}-R(1+\sum_{k=1}^{3}\beta_{k}cos(\eta_{k}(arctan(\frac{y}{x})-\theta_{k}))),$$
where the parameters are as follows:
\begin{equation}\label{A3}
\begin{split}
&\beta_{1}=0.3, \quad\;\ \, \eta_{1}=3, \quad\; \theta_{1}=0.5,\\
&\beta_{2}=-0.1, \quad \eta_{2}=4, \quad\, \theta_{2}=1.8,\\
&\beta_{3}=0.15, \quad\; \eta_{3}=7, \quad\, \theta_{3}=0, \quad R=0.483.
\end{split}
\end{equation}

\begin{figure}[H]
    \centering
    \includegraphics[width=16cm,height=11cm]{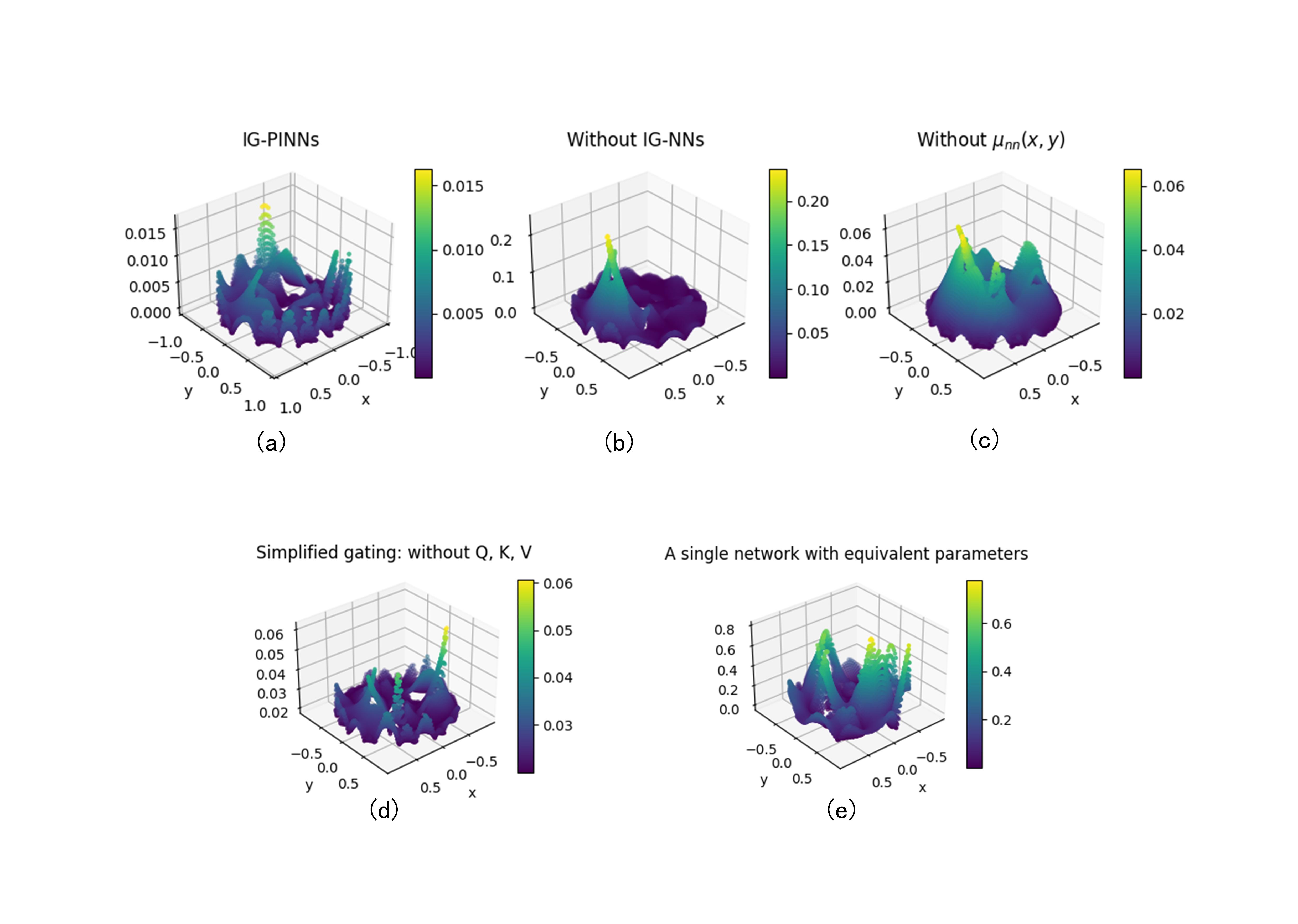}
   \caption{Ablation experiment, (\textbf{a}): Absolute errors of IG-PINNs. (\textbf{b}): Absolute errors without IG-NNs. (\textbf{c}): Absolute errors without $\mu_{nn}(x,y)$. (\textbf{d}): Absolute errors with simplified gating (without Q, K, V). (\textbf{e}): Absolute errors for a single network with equivalent parameters.}
    \label{AE_fig_error}
    \end{figure}

\begin{table}[H]
	\centering
	\caption{Comparison of errors in ablation experiment.}
	\label{AE_tab_error}
	\begin{tabular}{ccccc ccc}
		\hline\hline\noalign{\smallskip}	
    		& IG-PINNs & Without IG-NNs & Without $\mu_{nn}$ & Single network with & Simplified gating  \\
            & \ & \ & \ &  equivalent parameters & \ & \\
		\noalign{\smallskip}\hline\noalign{\smallskip}
		  $E_{M}$ & $1.14\times10^{-2}$ & $2.03\times10^{-1}$ & $5.97\times10^{-2}$ & $7.15\times10^{-1}$ & $6.12\times10^{-2}$\\
        $E_{L^2}$ & $9.86\times10^{-3}$ & $1.21\times10^{-1}$ & $7.02\times10^{-2}$ & $6.35\times10^{-1}$ & $5.73\times10^{-2}$\\
        $E_{jump}$ & $8.78\times10^{-3}$ & $1.01\times10^{-1}$ & $2.65\times10^{-2}$ & $1.01\times10^{0}$ & $3.04\times10^{-2}$\\
        $E_{flux}$ & $1.65\times10^{-3}$ & $6.73\times10^{-3}$ & $2.57\times10^{-3}$ &$9.72\times10^{-1}$ &$2.08\times10^{-3}$\\
		\noalign{\smallskip}\hline
	\end{tabular}
\end{table}

To evaluate the contribution of each component in the proposed model, we conducted systematic ablation studies, including:
\begin{itemize}
    \item \textbf{Full IG-PINNs},
    \item \textbf{IG-PINNs without IG-NNs, replacing IG-NNs with a fully connected neural network},
    \item \textbf{Direct approach with IG-NNs},
    \item \textbf{Single large network with equivalent parameters},
    \item \textbf{Simplified gating: \( Z^n = \sigma(W[H^{n-1};\phi(x)]) \) without \( Q, K, V \) machinery}.
\end{itemize}

We use 1,000 interior points, 1,200 outer boundary points, 120 internal boundary points, and 120 interface points as training data. The IG-PINNs architecture comprises three sub-networks: two IG-NNs and one fully connected neural network. The IG-NNs consist of a module with 32 neurons, and the fully connected neural network contains three hidden layers with 16 neurons. We set the maximum number of iterations to 50,000.

Figure \ref{AE_fig_error} shows the absolute errors obtained by IG-PINNs and other ablation experiments, with numerical results demonstrating that IG-PINNs provide the best approximation. Furthermore, Table \ref{AE_tab_error} offers a more detailed comparison. It can be observed that all error metrics of IG-PINNs are significantly lower than those obtained from other ablation experiments. The experimental results indicate that removing any network component leads to performance degradation in IG-PINNs. Across all error metrics, IG-PINNs clearly outperform the fully connected neural network with equivalent parameters. This shows that performance improvement mainly comes from the architectural design of the IG-PINNs, rather than merely from the increase in the parameter count.

\begin{table}[H]
	\centering
	\caption{Comparison of errors in Example \ref{E1}.}
	\label{E1_tab_error}
	\begin{tabular}{ccccc ccc}
		\hline\hline\noalign{\smallskip}	
		& PINNs & I-PINNs & M-PINNs & IG-PINNs  \\
		\noalign{\smallskip}\hline\noalign{\smallskip}
		  $E_{M}$ & $2.22\times10^{-1}$ & $5.31\times10^{-6}$ & $2.87\times10^{-6}$ & $1.18\times10^{-7}$ \\
          $E_{L^2}$ & $1.01\times10^{0}$ & $2.10\times10^{-5}$ & $1.05\times10^{-5}$ & $2.81\times10^{-7}$ \\
          Elapsed time & $21$(s) & $21$(s) & $26$(s) & $32$(s)\\
		\noalign{\smallskip}\hline
	\end{tabular}
\end{table}

\begin{figure}[H]
    \centering
    \includegraphics[width=14cm,height=9cm]{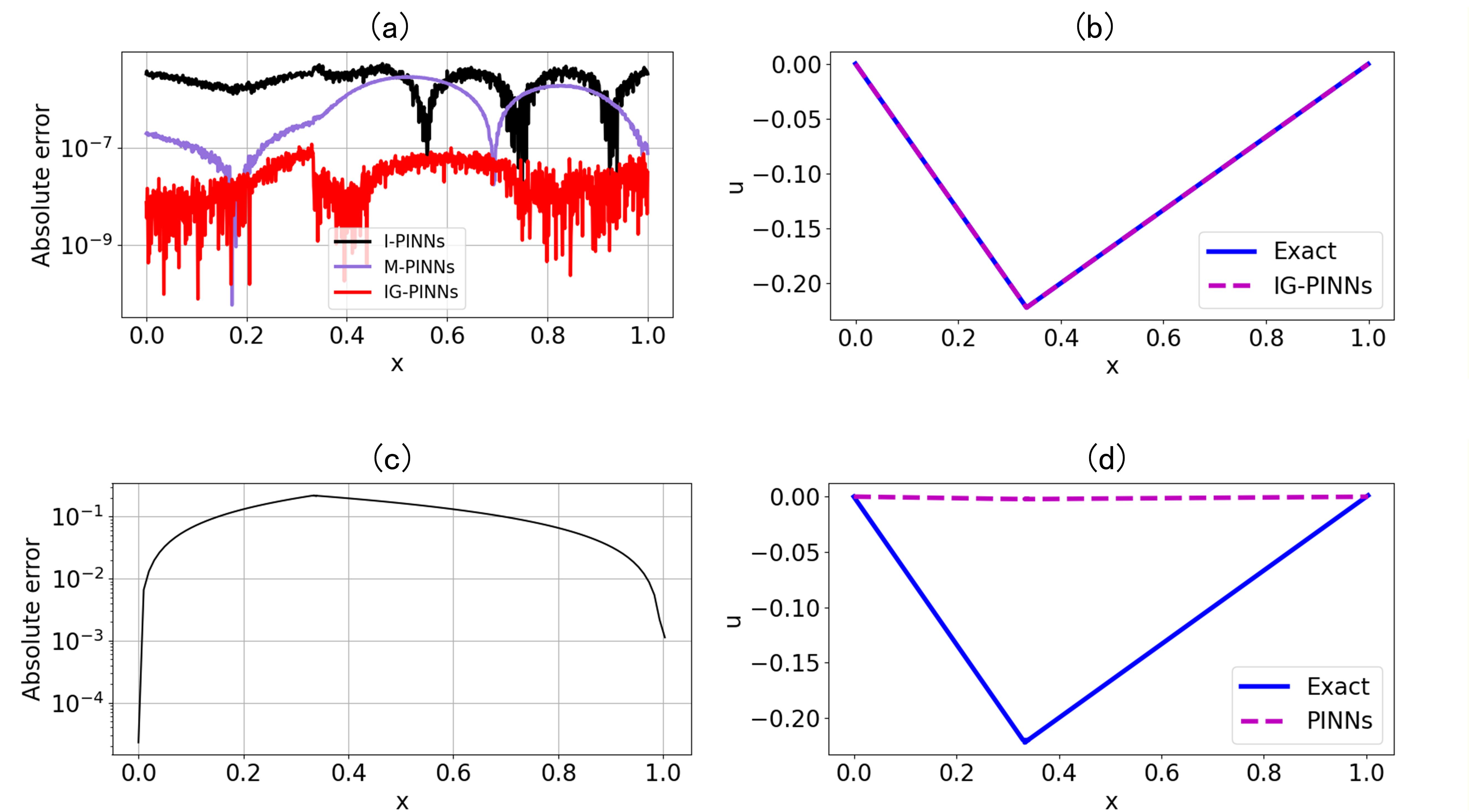}
    \caption{Example \ref{E1}, (\textbf{a}): Absolute errors of M-PINNs, I-PINNs and IG-PINNs. (\textbf{b}): IG-PINNs prediction, exact solution. (\textbf{c}): PINNs absolute errors. (\textbf{d}): PINNs prediction, exact solution.}
    \label{E1_fig_error}
    \end{figure}

\begin{figure}[H]
    \centering
    \includegraphics[width=10cm,height=6cm]{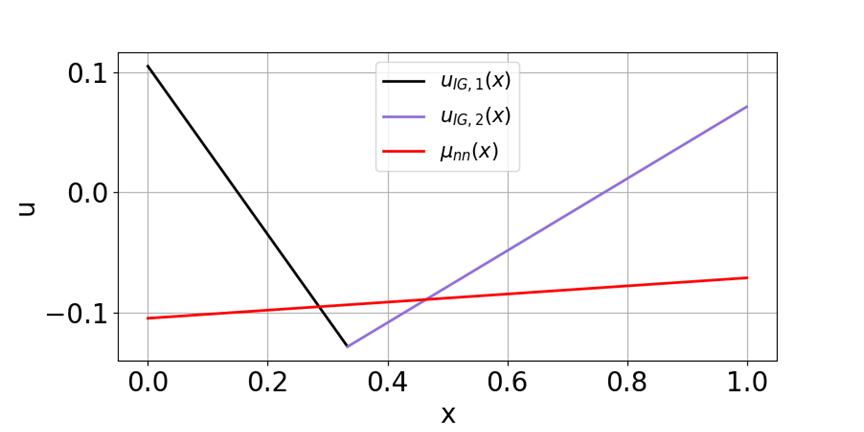}
    \caption{Example \ref{E1}, network components of IG-PINNs}
    \label{E1_fig_com}
    \end{figure}

\begin{example}\label{E1}
Consider the following $1$D Poisson equation:
\begin{equation}\label{3.12}
\left\{\begin{aligned}
\Delta u(x)&=0, \qquad x\in (0,1)/\{x_{\Gamma}\},\\
[\![u]\!](x_{\Gamma})&=0,\\
[\![\nabla u]\!](x_{\Gamma})&=1,\\
u(0)=u(1)&=0.\\
\end{aligned}\right.
\end{equation}
The exact solution:
\begin{equation}\label{3.13}
u(x)=\left\{
	\begin{aligned}
	(x_{\Gamma}-1)x, \qquad x\in [0,x_{\Gamma}],\\
	x_{\Gamma}(x-1), \qquad x\in [x_{\Gamma},1].\\
	\end{aligned}
	\right.
\end{equation}
\end{example}

In this problem, the solution is continuous throughout the domain, but its derivative is discontinuous at $x_{\Gamma}=\frac{1} {3}$. In order to demonstrate the effectiveness of the proposed algorithm, we choose I-PINNs, multi-domain physics-informed neural networks (M-PINNs) \cite{6} and PINNs as reference algorithms.

One module with $8$ neurons is used for IG-NNs. $\mu_{nn}(x)$ is a fully connected neural network with three hidden layers, each hidden layer containing 16 neurons. The level set function derived from the interface is set to $\phi(x)=x-\frac{1}{3}$. $100$ points in the interval $(0, 1)$, two boundary points and one interface point are used for training, and 1,000 test points are used to generate network prediction and calculate errors.  We set the maximum number of iterations to 15,000.

Figure \ref{E1_fig_error} (a) illustrates the absolute errors of IG-PINNs, I-PINNs and M-PINNs, where the maximum absolute errors of I-PINNs, M-PINNs and IG-PINNs are $5.31\times10^{-6}$, $2.87\times10^{ -6}$ and $1.18\times10^{-7}$, respectively. Figure \ref{E1_fig_com} shows the various network components of IG-PINNs. Figure \ref{E1_fig_error} (b) shows the IG-PINNs approximation and the exact solution, and it can be found that the IG-PINNs accurately approximate the exact solution. Figure \ref{E1_fig_error} (c), (d) show the absolute errors of the PINNs, the approximate solution, respectively. Since the activation function is inherently continuous, the PINNs fail to learn the discontinuity of the derivatives at the interface.  This results in a prediction of PINNs with large errors.

 Table \ref{E1_tab_error} shows the maximum absolute errors and the relative $L^{2}$ errors for IG-PINNs, I-PINNs and PINNs. By embedding interface information into neural networks, IG-PINNs achieve the best approximation compared with other algorithms. During the experiments, we found that the shallow neural network achieves a good approximation. Therefore, deeper networks are not used in the experiments.

\begin{figure}[H]
    \centering
    \includegraphics[width=11cm,height=8.5cm]{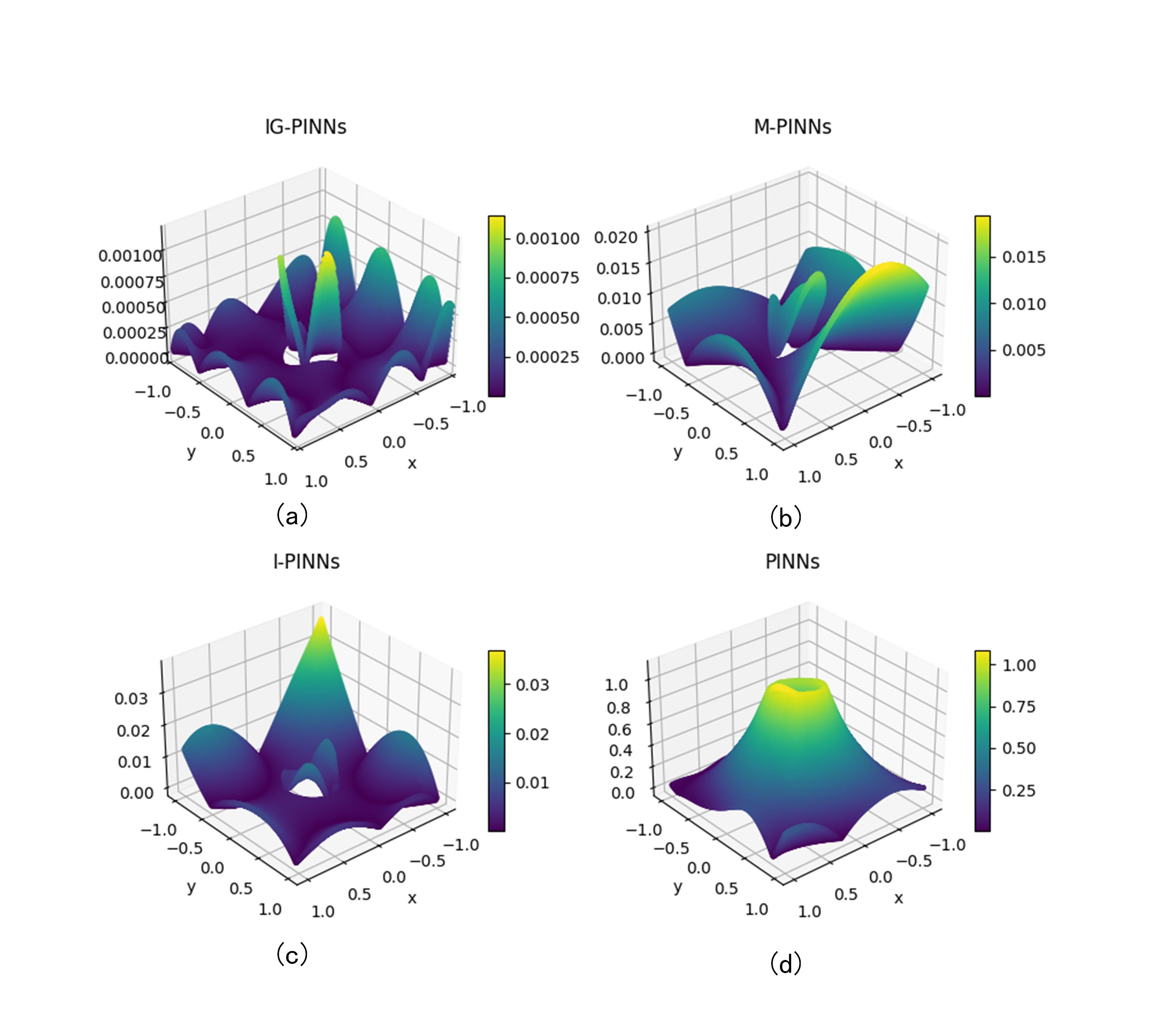}
    \caption{Example \ref{NN-FEM},  (\textbf{a}): IG-PINNs absolute errors. (\textbf{b}): M-PINNs absolute errors. (\textbf{c}): I-PINNs absolute errors. (\textbf{d}): PINNs absolute errors.}
    \label{NNFEM_fig_error}
 \end{figure}

 \begin{figure}[H]
    \centering
    \includegraphics[width=17cm,height=6cm]{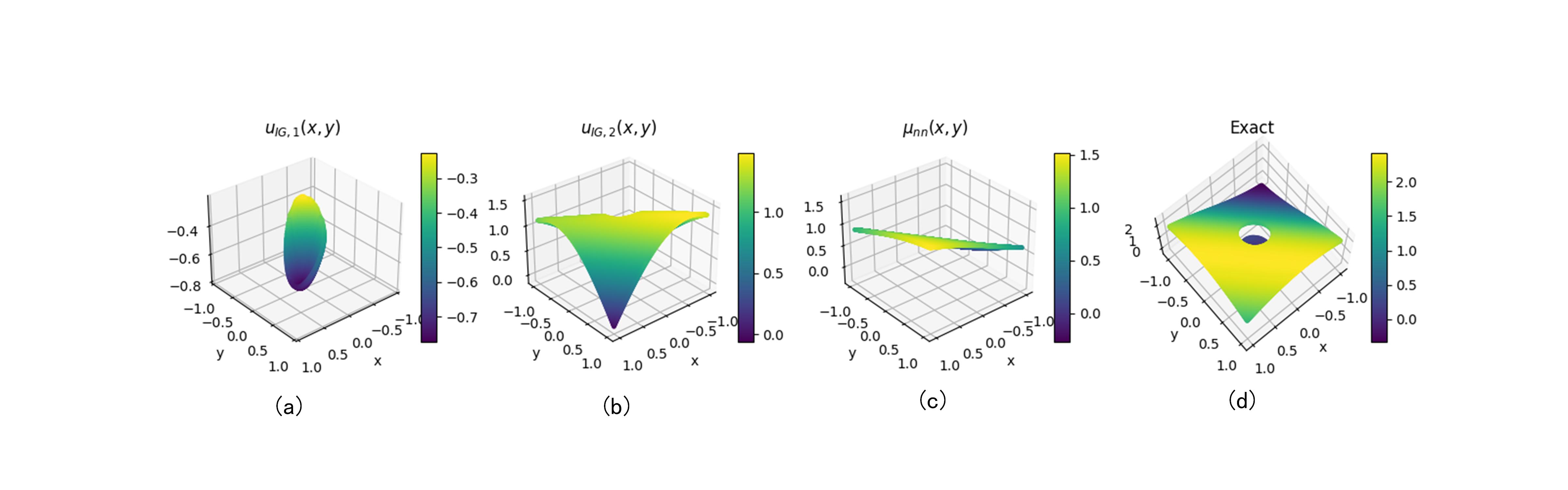}
    \caption{Example \ref{NN-FEM}, (\textbf{a}): Network component $u_{IG,1}(x,y)$. (\textbf{b}): Network component $u_{IG,2}(x,y)$. (\textbf{c}): Network component $\mu_{nn}(x,y)$. (\textbf{d}): Exact solution.}
    \label{NNFEM_fig_com}
 \end{figure}

\begin{example}\label{NN-FEM}
We consider a $2$D interface problem:
\begin{equation}\label{3.14}
\left\{\begin{aligned}
\Delta u(x,y) =f(x,y), &\qquad(x,y)\in \Omega/\Gamma,\\
u(x,y)=m(x,y), &\qquad(x,y)\in\partial\Omega,\\
[\![u]\!](x,y)=g(x,y), &\qquad(x,y)\in\Gamma,\\
([\![\nabla u]\!]\cdot\textbf{n})(x,y)=h(x,y),&\qquad(x,y)\in\Gamma,\\
\end{aligned}\right.
\end{equation}
where $\Omega=[-1,1]^{2}$ is divided into two subdomains $\Omega_{1},\Omega_{2}$ by the interface $\Gamma$:
$$\Gamma:x^2 + y^2=0.3^2.$$
The exact solution is defined as follows:
\begin{equation}\label{3.15}
u(x)=\left\{
	\begin{aligned}
	2x^{2}+3y^{2}, &\qquad (x,y)\in \Omega_{1},\\
	sin(x+y) +cos(x+y)+1,          &\qquad (x,y)\in \Omega_{2}.\\
	\end{aligned}
	\right.
\end{equation}
\end{example}

In this test, we used 500 interior points, 400 boundary points, and 2,000 interface points as training data. We employed a module with 32 neurons to construct IG-NNs.  The $\mu_{nn}(x,y)$ contains three hidden layers with 20 neurons. The maximum number of iterations is set to 20,000.

Figure \ref{NNFEM_fig_error} displays the absolute errors obtained by IG-PINNs and the reference algorithms. Figure \ref{NNFEM_fig_com} shows the network components of IG-PINNs. Table \ref{NNFEN_tab_error} provides a comparison of various errors. It can be observed that IG-PINNs achieve the best approximation among all algorithms. Table \ref{NNFEM_tab_module} reports the relative $L^{2}$ errors for different numbers of modules and neurons. The numerical results indicate that IG-NNs with one module and shallow fully connected neural networks can provide a good approximation. Table \ref{NNFEM_tab_inter} reports the relative $L^{2}$ errors for different numbers of interface points. As the number of interface points increases, the approximation of IG-PINNs is improved.

\begin{figure}[H]
    \centering
    \includegraphics[width=14cm,height=6cm]{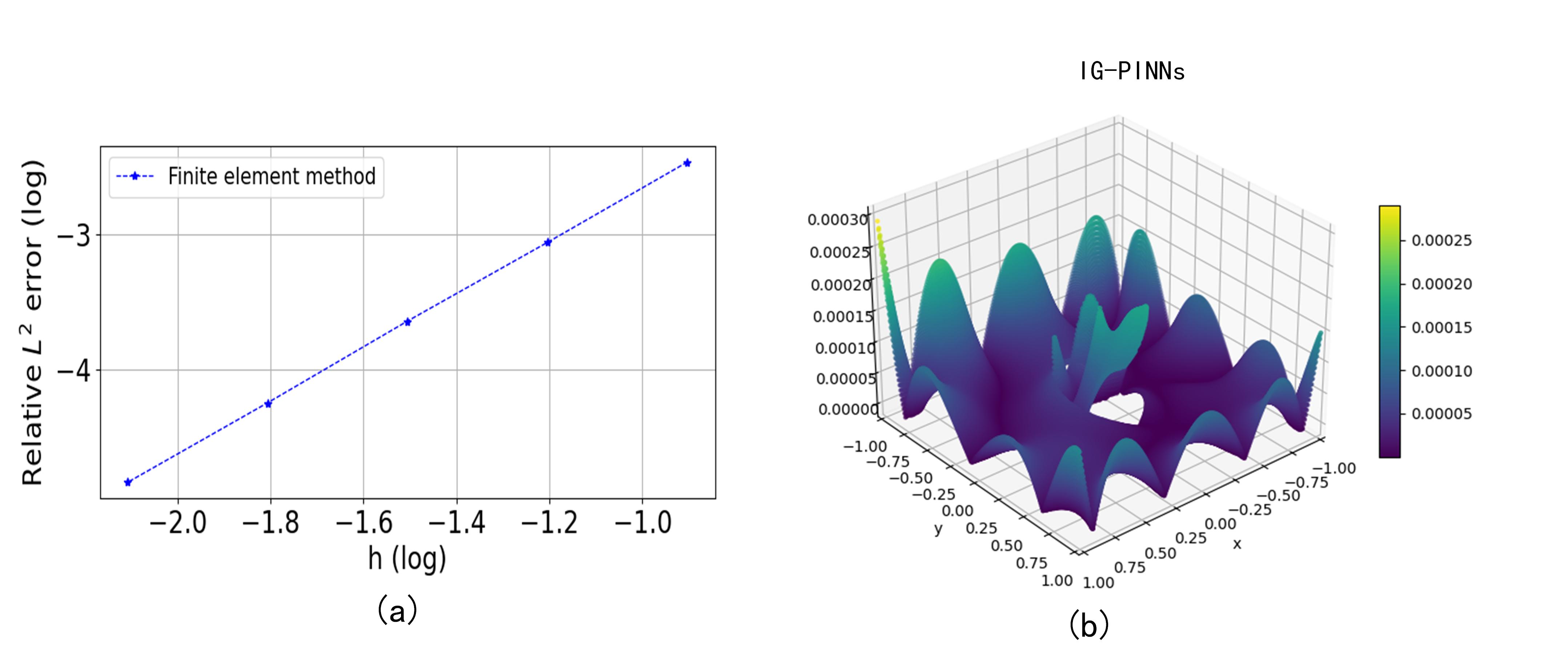}
    \caption{Example \ref{NN-FEM}, (a): The relative $L^{2}$ errors obtained by FEM \cite{57}, where $h$ represents the step size. (b): IG-PINNs absolute errors.}
    \label{NNFEM_FEM}
 \end{figure}

Many traditional numerical algorithms have made significant progress in solving interface problems, including finite volume method \cite{11}, immersed boundary method \cite{12} and ghost fluid method \cite{13,14}. To further evaluate the performance of IG-PINNs in solving interface problems, the finite element method (FEM) \cite{57} is used as a reference benchmark in this case. Figure \ref{NNFEM_FEM} (a) displays the relative $L^{2}$ errors obtained by the FEM. As the mesh is refined, the relative $L^{2}$ error of FEM decreases to $1.45\times10^{-5}$, whereas that of IG-PINNs is $1.01\times10^{-4}$. To further improve the accuracy of IG-PINNs, we introduced a learning rate decay mechanism. The initial learning rate is set to 0.001, and the decay steps and decay rate are set to 1,000 and 0.95, respectively. This strategy allows the model to converge quickly with a larger learning rate in the early stages of training, followed by a gradual reduction in the learning rate to fine-tune the parameters, thereby avoiding oscillations and enhancing the model's generalization ability. The experimental results are presented in Figure \ref{NNFEM_FEM} (b), where the relative $L^{2}$ error is $2.74\times10^{-5}$. It can be observed that there is still a significant gap between IG-PINNs and the finite element method. However, FEM often encounters difficulties in dealing with interface conditions involving complex geometric structures and dynamic boundaries, and it is difficult to solve high-dimensional problems. Neural networks hold advantages in addressing these challenges.

\begin{table}[H]
	\centering
	\caption{Comparison of errors in Example \ref{NN-FEM}.}
	\label{NNFEN_tab_error}
	\begin{tabular}{ccccc ccc}
		\hline\hline\noalign{\smallskip}	
    		& PINNs & I-PINNs & M-PINNs & IG-PINNs  \\
		\noalign{\smallskip}\hline\noalign{\smallskip}
		  $E_{M}$ & $1.08\times10^{0}$ & $3.11\times10^{-2}$ & $1.53\times10^{-2}$ & $1.04\times10^{-3}$\\
        $E_{L^2}$ & $2.84\times10^{-1}$ & $4.27\times10^{-3}$ & $3.04\times10^{-3}$ & $1.01\times10^{-4}$ \\
        $E_{jump}$ & $9.24\times10^{-1}$ & $3.65\times10^{-3}$ & $8.48\times10^{-3}$ & $7.14\times10^{-5}$ \\
        $E_{flux}$ & $7.41\times10^{-1}$ & $5.17\times10^{-3}$ & $1.72\times10^{-3}$ & $2.24\times10^{-4}$ \\
        Elapsed time &$209$(s) & $203$(s) & $256$(s) & $268$(s) \\
		\noalign{\smallskip}\hline
	\end{tabular}
\end{table}

\begin{table}[H]
	\centering
	\caption{Comparison of relative $L^{2}$ errors in Example \ref{NN-FEM}.}
	\label{NNFEM_tab_module}
	\begin{tabular}{ccccc ccc}
		\hline\hline\noalign{\smallskip}	
    	Neurons & \ & 8 & 16 & 32  & 48\\
		\noalign{\smallskip}\hline\noalign{\smallskip}
        &$1$ & $1.92\times10^{-3}$ & $5.71\times10^{-4}$ & $1.01\times10^{-4}$ & $1.27\times10^{-4}$ \\
        Modules &$2$ & $7.13\times10^{-4}$ & $4.26\times10^{-4}$ & $1.83\times10^{-4}$ & $2.17\times10^{-4}$ \\
        &$3$ & $7.71\times10^{-4}$ & $9.61\times10^{-5}$ & $1.07\times10^{-4}$ & $1.12\times10^{-4}$ \\
		\noalign{\smallskip}\hline
	\end{tabular}
\end{table}

\begin{table}[H]
	\centering
	\caption{Errors for different numbers of interface points in Example \ref{NN-FEM}.}
	\label{NNFEM_tab_inter}
	\begin{tabular}{ccccc ccc}
		\hline\hline\noalign{\smallskip}	
    	Interface points & 500 & 1000 & 2000  & 3000\\
		\noalign{\smallskip}\hline\noalign{\smallskip}
        $E_{L^2}$ & $3.16\times10^{-3}$ & $6.45\times10^{-4}$ & $1.01\times10^{-4}$ & $1.42\times10^{-4}$ \\
        $E_{jump}$ & $4.17\times10^{-3}$ & $6.27\times10^{-4}$ & $7.14\times10^{-5}$ & $8.26\times10^{-5}$ \\
        $E_{flux}$ & $8.37\times10^{-3}$ & $9.17\times10^{-4}$ & $2.24\times10^{-4}$ & $3.13\times10^{-4}$ \\
		\noalign{\smallskip}\hline
	\end{tabular}
\end{table}

\begin{table}[H]
	\centering
	\caption{Comparison of errors in Example \ref{E2}.}
	\label{E2_tab_error}
	\begin{tabular}{ccccc ccc}
		\hline\hline\noalign{\smallskip}	
    		& PINNs & I-PINNs & M-PINNs & IG-PINNs  \\
		\noalign{\smallskip}\hline\noalign{\smallskip}
		  $E_{M}$ & $1.52\times10^{0}$ & $6.17\times10^{-1}$ & $4.41\times10^{-1}$ & $1.22\times10^{-2}$\\
        $E_{L^2}$ & $1.71\times10^{-1}$ & $2.69\times10^{-2}$ & $1.07\times10^{-2}$ & $4.08\times10^{-4}$ \\
        $E_{jump}$ & $5.84\times10^{-2}$ & $1.64\times10^{-2}$ & $8.05\times10^{-3}$ & $4.49\times10^{-4}$ \\
        $E_{flux}$ & $1.54\times10^{-2}$ & $6.27\times10^{-3}$ & $8.12\times10^{-4}$ & $1.52\times10^{-4}$ \\
        Elapsed time &$307$(s) & $305$(s) & $403$(s) & $353$(s) \\
		\noalign{\smallskip}\hline
	\end{tabular}
\end{table}

\begin{table}[H]
	\centering
	\caption{Errors for different numbers of interface points in Example \ref{E2}.}
	\label{E2_tab_inter}
	\begin{tabular}{ccccc ccc}
		\hline\hline\noalign{\smallskip}	
    	Interface points & 200 & 400 & 600  & 800\\
		\noalign{\smallskip}\hline\noalign{\smallskip}
        $E_{L^2}$ & $1.18\times10^{-3}$ & $7.12\times10^{-4}$ & $4.08\times10^{-4}$ & $4.54\times10^{-4}$ \\
        $E_{jump}$ & $6.86\times10^{-3}$ & $4.59\times10^{-4}$ & $4.49\times10^{-4}$ & $5.06\times10^{-4}$ \\
        $E_{flux}$ & $6.28\times10^{-4}$ & $2.19\times10^{-4}$ & $1.52\times10^{-4}$ & $2.23\times10^{-4}$ \\
		\noalign{\smallskip}\hline
	\end{tabular}
\end{table}

\begin{table}[H]
	\centering
	\caption{Comparison of relative $L^{2}$ errors in Example \ref{E2}.}
	\label{E2_tab_module}
	\begin{tabular}{ccccc ccc}
		\hline\hline\noalign{\smallskip}	
    	Neurons & \ & 8 & 16 & 32  & 48\\
		\noalign{\smallskip}\hline\noalign{\smallskip}
        &$1$ & $2.71\times10^{-2}$ & $2.42\times10^{-3}$ & $4.08\times10^{-4}$ & $4.59\times10^{-4}$ \\
        Modules &$2$ & $9.13\times10^{-3}$ & $7.31\times10^{-4}$ & $4.47\times10^{-4}$ & $4.32\times10^{-4}$ \\
        &$3$ & $5.49\times10^{-3}$ & $6.56\times10^{-4}$ & $4.74\times10^{-4}$ & $5.87\times10^{-4}$ \\
		\noalign{\smallskip}\hline
	\end{tabular}
\end{table}

\begin{example}\label{E2}
To evaluate the performance of IG-PINNs in interface problems with non-smooth interfaces, we consider a 2D elliptic interface problem with corners:
\begin{equation}\label{3.14}
\left\{\begin{aligned}
\nabla\cdot(\alpha(x,y)(\nabla u(x,y)))=f(x,y), &\qquad(x,y)\in \Omega/\Gamma,\\
u(x,y)=m(x,y), &\qquad(x,y)\in\partial\Omega,\\
[\![u]\!](x,y)=g(x,y), &\qquad(x,y)\in\Gamma,\\
([\![\alpha\nabla u]\!]\cdot\textbf{n})(x,y)=h(x,y),&\qquad(x,y)\in\Gamma,\\
\end{aligned}\right.
\end{equation}
where $\Omega=[-1,1]^{2}$ is divided into two subdomains $\Omega_{1},\Omega_{2}$ by the interface $\Gamma$:
$$\Gamma:max(|x|,|y|)=0.5.$$
The exact solution is defined as follows:
\begin{equation}\label{3.15}
u(x)=\left\{
	\begin{aligned}
	10e^{x^{2}+y^{2}}, &\qquad (x,y)\in \Omega_{1},\\
	3(x+y),          &\qquad (x,y)\in \Omega_{2},\\
	\end{aligned}
	\right.
\end{equation}
and
\begin{equation}\label{E2A}
\begin{split}
\alpha(x,y)=\left\{
	\begin{aligned}
	10^{-2}, &\qquad(x,y)\in \Omega_{1},\\
	10^{2}, &\qquad(x,y)\in \Omega_{2}.\\
	\end{aligned}
	\right.
\end{split}
\end{equation}
\end{example}

In this case, the level set function derived from the interface is set as follows:
$$\phi(x,y)=(x^{2}-0.5^{2})(y^{2}-0.5^{2}).$$
For network training, 1,000 interior points, 400 boundary points, and 600 interface points are sampled by the Latin hypercube sampling (LHS) \cite{53}. We use 40,000 test points to generate neural network predictions and calculate errors for reference algorithms. During the experiments, two modules are used to construct the IG-NNs. $\mu_{nn}(x,y)$ is a fully connected neural network with four hidden layers, each hidden layer containing 32 neurons.  We set the maximum number of iterations to 20,000.

\begin{figure}[H]
    \centering
    \includegraphics[width=13cm,height=12cm]{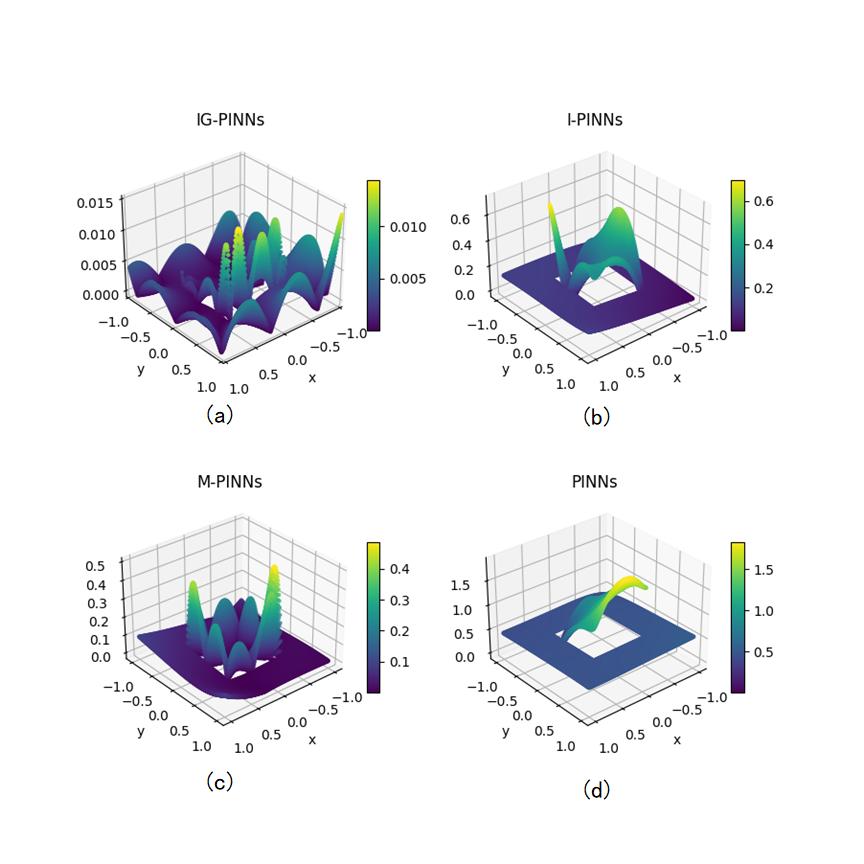}
    \caption{Example \ref{E2},  (\textbf{a}): IG-PINNs absolute errors. (\textbf{b}): I-PINNs absolute errors. (\textbf{c}): M-PINNs absolute errors. (\textbf{d}): PINNs absolute errors.}
    \label{E2_fig_error}
 \end{figure}

\begin{figure}[H]
    \centering
    \includegraphics[width=17cm,height=6cm]{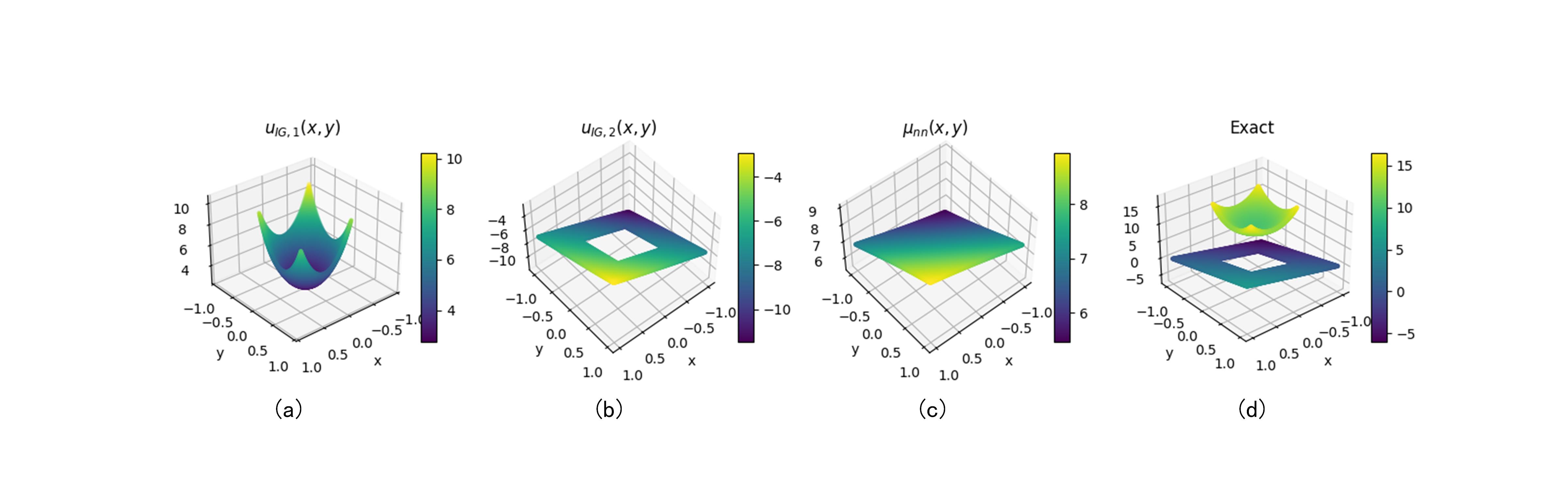}
    \caption{Example \ref{E2}, (\textbf{a}): Network component $u_{IG,1}(x,y)$. (\textbf{b}): Network component $u_{IG,2}(x,y)$. (\textbf{c}): Network component $\mu_{nn}(x,y)$. (\textbf{d}): Exact solution.}
    \label{E2_fig_com}
 \end{figure}

Figure \ref{E2_fig_error} shows the absolute errors of the four neural network algorithms. The results show that IG-PINNs accurately approximate the solution of problem  compared with other algorithms. Figure \ref{E2_fig_com} shows the various network components of IG-PINNs. Table \ref{E2_tab_error} shows the various errors obtained by the four neural network algorithms. Although we observed that I-PINNs had the shortest training time, they performed poorly in solving interface problems with corners (\ref{3.14}), exhibiting a relative $L^{2}$ error of $2.69\times10^{-2}$. In contrast, IG-PINNs achieved better approximation results, with a relative $L^{2}$ error of only $4.08\times10^{-4}$, which is significantly lower than that of the former. In subdomain $\Omega_{1}$, the solution is only constrained by the interface conditions and multiple components in the loss function interact. A good neural network representation is obtained depending on the interface conditions being sufficiently optimized. By embedding interface information into the network structure, IG-PINNs achieve a better approximation compared with PINNs, I-PINNs and M-PINNs.

We test the performance of IG-PINNs at different numbers of modules and training points. Table \ref{E2_tab_module} provides a detailed comparison. Numerical results indicate that IG-PINNs with a module can achieve a good approximation. Furthermore, we tested the model performance with different numbers of interface points, and the numerical results are presented in Table \ref{E2_tab_inter}.

\begin{example}\label{E3}
In this example, considering a $2$D interface problem with a star-shaped interface:
\begin{equation}\label{3.16}
\left\{\begin{aligned}
\Delta u(x,y)=f(x,y), &\qquad(x,y)\in \Omega/ \Gamma,\\
[\![u]\!](x,y)=g(x,y), &\qquad(x,y)\in\Gamma,\\
u(x,y)=m(x,y), &\qquad(x,y)\in\partial\Omega,\\
([\![\nabla u]\!]\cdot\textbf{n})(x,y)=h(x,y), &\qquad(x,y)\in\Gamma.\\
\end{aligned}\right.
\end{equation}
The interface is described by a star curve:
\begin{equation}\label{3.17}
\begin{split}
\Gamma=\{(x,y):\phi(x,y)=0,\phi(x,y)=\sqrt{(x^{2}+y^{2})}-(0.5+sin(arctanh(\theta))/7\}.
\end{split}
\end{equation}
The exact solution:
\begin{equation}\label{3.19}
\begin{split}
u(x,y)=\left\{
	\begin{aligned}
	x^{2}+y^{2}, &\qquad(x,y)\in \Omega_{1},\\
	0.1(x^{2}+y^{2})^{2}-0.01log(2\sqrt{x^{2}+y^{2}}), &\qquad(x,y)\in \Omega_{2}.\\
	\end{aligned}
	\right.
\end{split}
\end{equation}
\end{example}

In this test, we use $400$ boundary points, $500$ interior points and $600$ interface points as training points. The number of iterations reached $10^{5}$. We use two modules with $48$ neurons to build IG-NNs. $\mu_{nn}(x,y)$ is a fully connected neural network with three hidden layers, each hidden layer containing 20 neurons.

\begin{figure}[H]
    \centering
    \includegraphics[width=13cm,height=11cm]{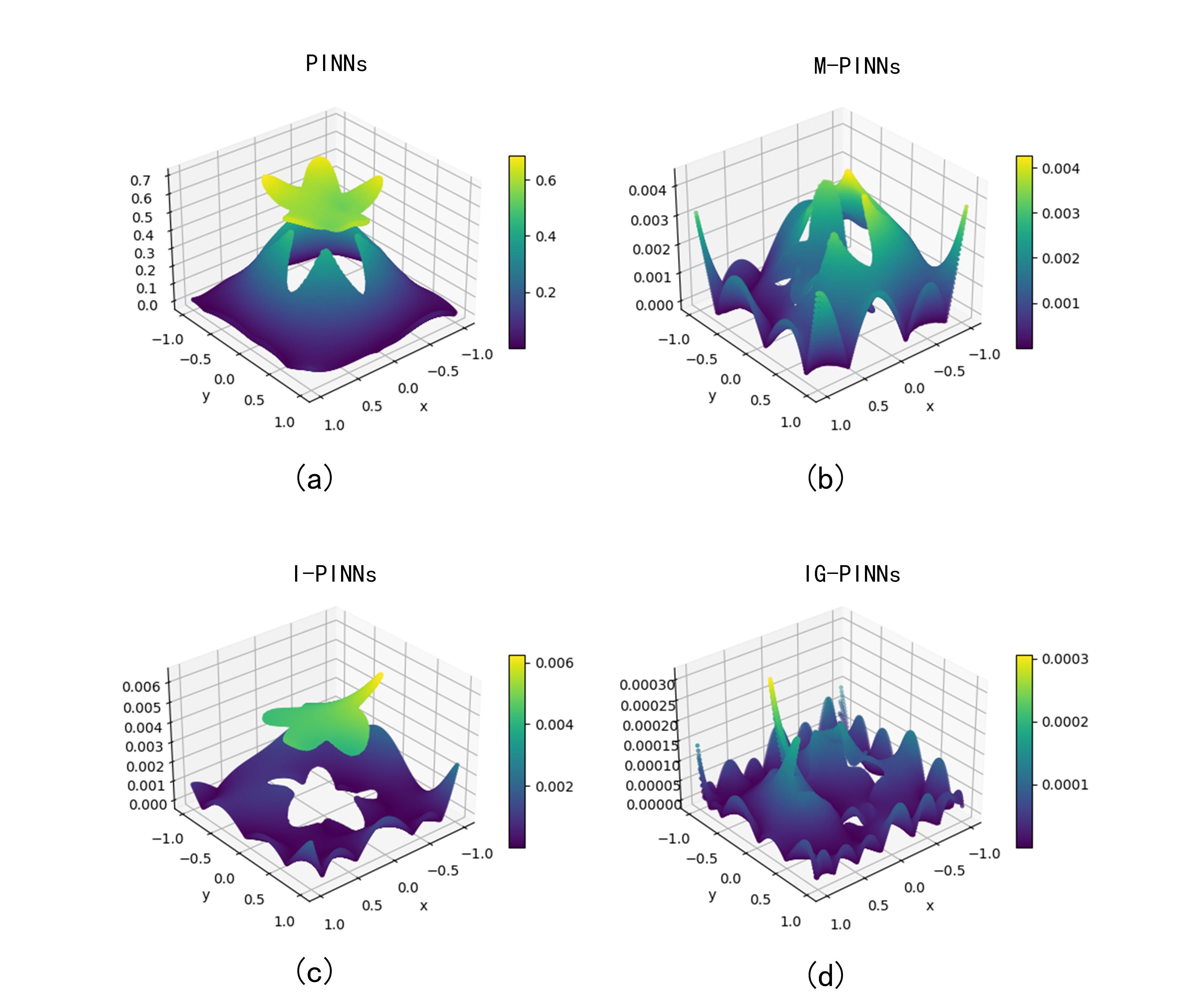}
    \caption{Example \ref{E3},  (\textbf{a}): PINNs absolute errors. (\textbf{b}): M-PINNs absolute errors. (\textbf{c}): I-PINNs absolute errors. (\textbf{d}): IG-PINNs absolute errors.}
    \label{E3_fig_error}
    \end{figure}

\begin{figure}[H]
    \centering
    \includegraphics[width=17cm,height=6cm]{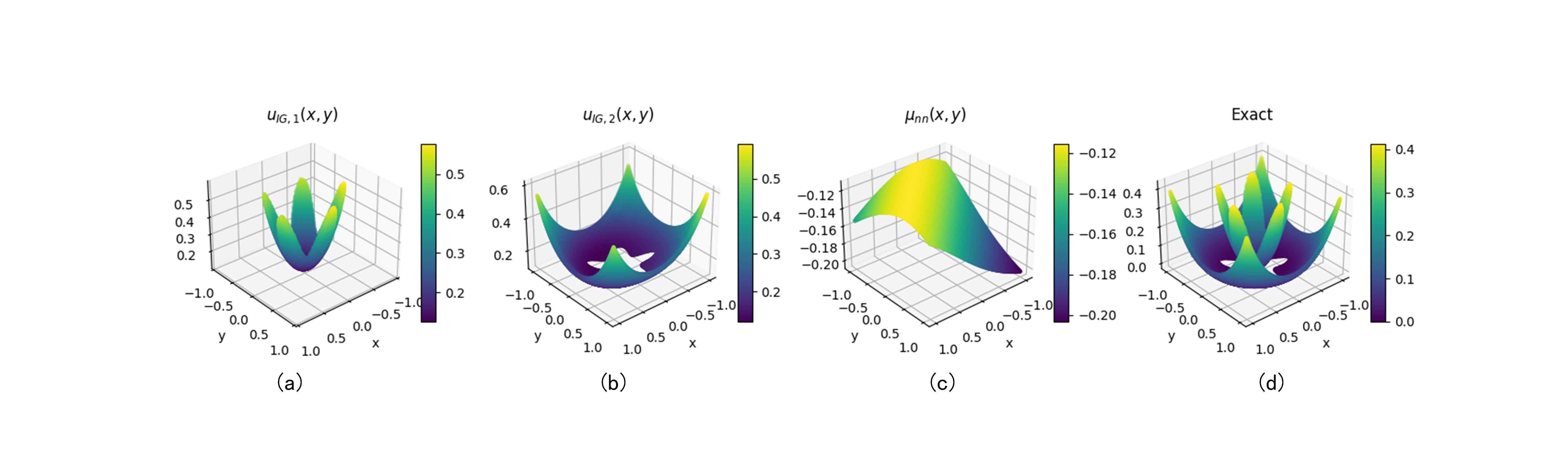}
    \caption{Example \ref{E3}, (\textbf{a}): Network component $u_{IG,1}(x,y)$. (\textbf{b}): Network component $u_{IG,2}(x,y)$. (\textbf{c}): Network component $\mu_{nn}(x,y)$. (\textbf{d}): Exact solution.}
    \label{E3_fig_com}
    \end{figure}

Figure \ref{E3_fig_error} shows the absolute errors obtained by different neural network algorithms. Figure \ref{E3_fig_com} shows the various network components of IG-PINNs. Note that the solution jumps across the interface, and the solution in subdomain $\Omega_{1}$ is constrained only by the interface conditions. A single neural network with continuous activations function fails to capture discontinuities across the interface, which results in a PINNs approximation with large errors. Although the approximations of I-PINNs and M-PINNs are improved compared with PINNs, the interface conditions are not properly optimized in I-PINNs and M-PINNs, resulting in error propagation to subdomain $\Omega_{1}$. IG-PINNs with interface information provide a more accurate approximation compared with I-PINNs and M-PINNs.

Table \ref{E3_tab_error} provides a detailed comparison of errors. The maximum absolute errors of  IG-PINNs, I-PINNs and M-PINNs are $3.05 \times10^{-4}$, $6.51 \times10^{-3}$ and $4.26 \times10^{-3}$ respectively. The relative $L^{2}$ error of IG-PINNs is an order of magnitude lower than that of I-PINNs and M-PINNs. It should be noted that IG-PINNs require more training time compared to the reference methods. We test the performance of IG-PINNs at different numbers of modules and training points. Table \ref{E3_tab_module} provides a detailed comparison. Furthermore, Table \ref{E3_tab_inter} provides a detailed error comparison for different numbers of interface points. It should be noted that although IG-PINNs achieve a good approximation, they require more training time compared to other algorithms.


\begin{table}[H]
	\centering
	\caption{Comparison of errors in Example \ref{E3}.}
	\label{E3_tab_error}
	\begin{tabular}{ccccc}
		\hline\hline\noalign{\smallskip}	
		& PINNs & I-PINNs & M-PINNs & IG-PINNs  \\
		\noalign{\smallskip}\hline\noalign{\smallskip}
		  $E_{M}$ & $6.85\times10^{-1}$ & $6.51\times10^{-3}$ & $4.26\times10^{-3}$ & $3.05\times10^{-4}$ \\
          $E_{L^2}$ & $2.48\times10^{0}$ & $2.23\times10^{-2}$ & $1.34\times10^{-2}$ & $4.28\times10^{-4}$ \\
          $E_{jump}$ & $7.46\times10^{-2}$ & $6.37\times10^{-3}$ & $4.56\times10^{-3}$ & $8.35\times10^{-4}$ \\
          $E_{flux}$ & $5.17\times10^{-2}$ & $5.59\times10^{-3}$ & $3.61\times10^{-3}$ & $7.68\times10^{-4}$ \\
          Elapsed time & $486$(s) & $493$(s) & $544$(s) & $605$(s) \\
		\noalign{\smallskip}\hline
	\end{tabular}
\end{table}

\begin{table}[H]
	\centering
	\caption{Errors for different numbers of interface points in Example \ref{E3}.}
	\label{E3_tab_inter}
	\begin{tabular}{ccccc ccc}
		\hline\hline\noalign{\smallskip}	
    	Interface points & 200 & 400 & 600  & 800\\
		\noalign{\smallskip}\hline\noalign{\smallskip}
        $E_{L^2}$ & $5.71\times10^{-2}$ & $3.81\times10^{-3}$ & $4.28\times10^{-4}$ & $4.87\times10^{-4}$ \\
        $E_{jump}$ & $9.15\times10^{-3}$ & $3.94\times10^{-3}$ & $8.35\times10^{-4}$ & $9.27\times10^{-4}$ \\
        $E_{flux}$ & $1.62\times10^{-2}$ & $2.09\times10^{-3}$ & $7.68\times10^{-4}$ & $8.73\times10^{-4}$ \\
		\noalign{\smallskip}\hline
	\end{tabular}
\end{table}

\begin{table}[H]
	\centering
	\caption{Relative $L^{2}$ errors for different numbers of modules and training points in Example \ref{E3}.}
	\label{E3_tab_module}
	\begin{tabular}{ccccc ccc}
		\hline\hline\noalign{\smallskip}	
    	Neurons & \ & 8 & 16 & 32  & 48\\
		\noalign{\smallskip}\hline\noalign{\smallskip}
        &$1$ & $2.51\times10^{-1}$ & $6.72\times10^{-2}$ & $2.48\times10^{-2}$ & $1.57\times10^{-2}$ \\
        Modules &$2$ & $2.73\times10^{-2}$ & $3.49\times10^{-3}$ & $4.28\times10^{-4}$ & $5.31\times10^{-4}$ \\
        &$3$ & $3.18\times10^{-2}$ & $1.43\times10^{-3}$ & $6.31\times10^{-4}$ & $5.22\times10^{-4}$ \\
		\noalign{\smallskip}\hline
	\end{tabular}
\end{table}

\begin{example}\label{E4}
We consider a $3$D elliptic interface problem with discontinuous coefficients in the computational domain $\Omega=[-1,1]^{3}$.
\begin{equation}\label{3.21}
\left\{\begin{aligned}
\nabla\cdot(\alpha(x,y,z)(\nabla u(x,y,z)))=f(x,y,z), &\qquad(x,y,z)\in \Omega/\Gamma,\\
([\![\alpha\nabla u]\!]\cdot\textbf{n})(x,y,z)=h(x,y,z), &\qquad(x,y,z)\in\Gamma,\\
[\![u]\!](x,y,z)=g(x,y,z), &\qquad(x,y,z)\in\Gamma,\\
u(x,y,z)=m(x,y,z), &\qquad(x,y,z)\in \partial\Omega.\\
\end{aligned}\right.
\end{equation}
An ellipsoidal interface $\Gamma=\{(x,y,z):\phi(x,y,z)=0,\phi(x,y,z)=2x^{2}+3y^{2}+6z^{2}-1.69\}$ divides the domain into non-overlapping subdomains $\Omega_{1}=\{(x,y,z):\phi(x,y,z)\leq0\}$, $\Omega_{2}=\{(x,y,z):\phi(x,y,z)>0\}$. The exact solution is piecewise defined:
\begin{equation}\label{3.22}
u(x,y,z)=\left\{
	\begin{aligned}
	& sin(2x)cos(2y)e^{z}, \quad\quad\quad\quad\quad\quad\quad\quad (x,y,z)\in \Omega_{1},\\
	&(16(\frac{y-x} {3})^{5}-20(\frac{y-x} {3})^{3}\\
    & +5(\frac{y-x} {3})) log(x+y+3)cos(z), \,\quad (x,y,z)\in \Omega_{2}.\\
	\end{aligned}
	\right.
\end{equation}
The coefficient $\alpha(x,y,z)$ jumps across the interface $\Gamma$:
\begin{equation}\label{3.23}
\alpha=\left\{
	\begin{aligned}
	10(1+0.2cos(2\pi(x+y))sin(2\pi(x-y))cos(z)), &\qquad (x,y,z)\in \Omega_{1},\\
	1, &\qquad(x,y,z)\in \Omega_{2}.\\
	\end{aligned}
	\right
	.
\end{equation}
\end{example}

In this test, we use LHS to generate $1500$ training points, including $500$ interior points, $600$ boundary points and $400$ interface points. To demonstrate the effectiveness of the proposed method, we use I-PINNs and M-PINNs as reference algorithms. In IG-PINNs, we use a module with $32$ neurons to construct the IG-NNs, and the level set function $\phi(x,y,z)=2x^{2}+3y^{2}+6z^{2}-1.69$. We use a fully connected neural network with three hidden layers as $\mu_{nn}(x,y,z)$. The number of iterations during the training process is set to $100,000$.  

\begin{figure}[H]
    \centering
    \includegraphics[width=14cm,height=7cm]{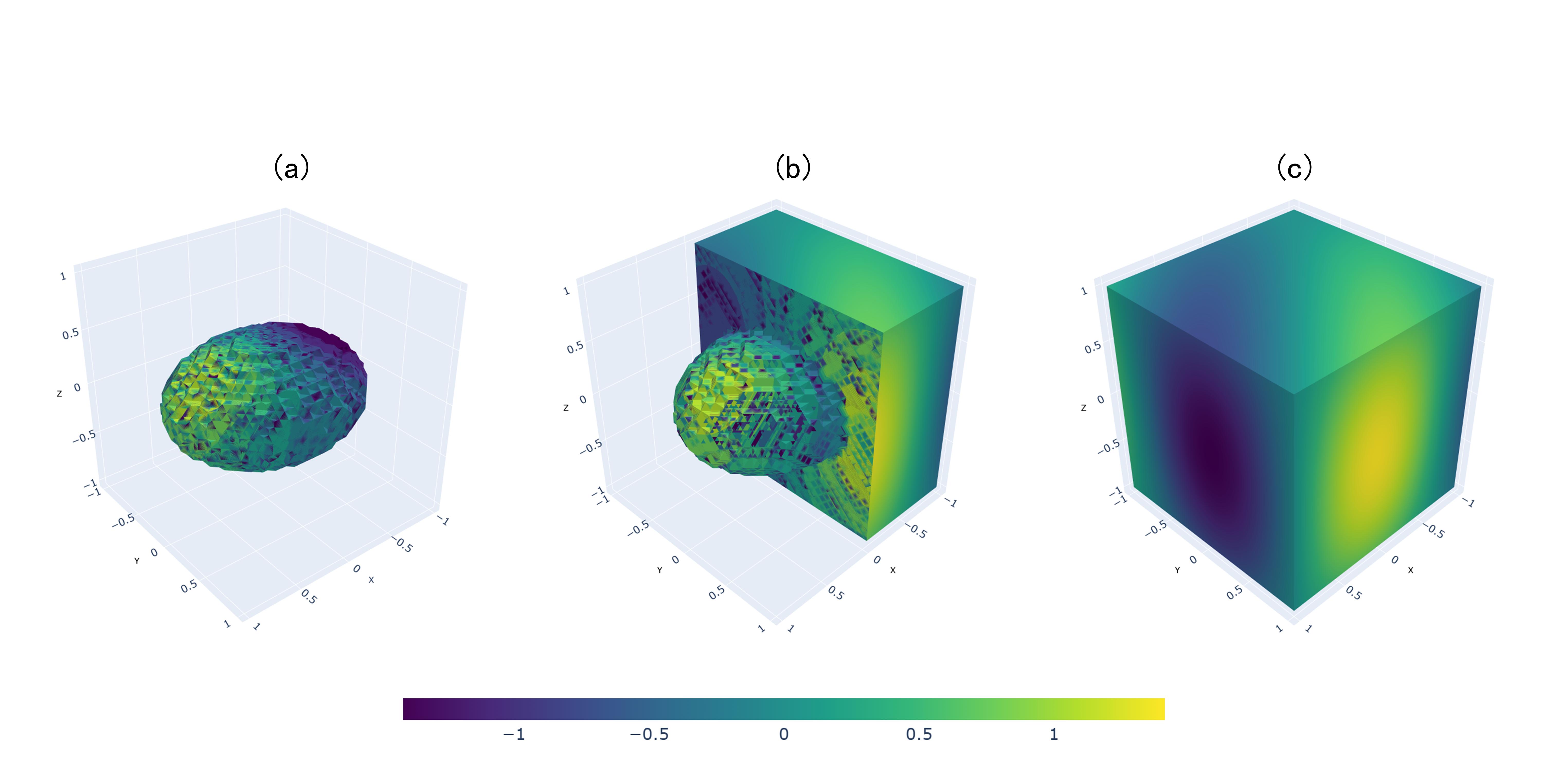}
    \caption{Example \ref{E4}, (\textbf{a}): The exact solution at the interface. (\textbf{b}): The exact solution inside the computational domain. (\textbf{c}): The exact solution at the boundary of the computational domain.}
    \label{E4_fig_solution}
\end{figure}

\begin{table}[H]
	\centering
	\caption{Comparison of errors in Example \ref{E4}.}
	\label{E4_tab_error}
	\begin{tabular}{ccccc ccc}
		\hline\hline\noalign{\smallskip}	
    		& I-PINNs & M-PINNs & IG-PINNs  \\
		\noalign{\smallskip}\hline\noalign{\smallskip}
		  $E_{M}$ & $2.07\times10^{-1}$ & $8.61\times10^{-2}$ & $1.75\times10^{-2}$ \\
        $E_{L^2}$ & $8.91\times10^{-2}$ & $5.23\times10^{-2}$ & $4.11\times10^{-3}$ \\
        $E_{jump}$ & $3.74\times10^{-2}$ & $1.53\times10^{-2}$ & $3.35\times10^{-3}$ \\
        $E_{flux}$ & $2.85\times10^{-2}$ & $6.47\times10^{-4}$ & $4.49\times10^{-4}$ \\
        Elapsed time & $835$(s) & $813$(s) & $1424$(s) \\
		\noalign{\smallskip}\hline
	\end{tabular}
\end{table}

\begin{table}[H]
	\centering
	\caption{Errors for different numbers of interface points in Example \ref{E4}.}
	\label{E4_tab_inter}
	\begin{tabular}{ccccc ccc}
		\hline\hline\noalign{\smallskip}	
    	Interface points & 100 & 400 & 900  & 1200\\
		\noalign{\smallskip}\hline\noalign{\smallskip}
        $E_{L^2}$ & $1.73\times10^{-2}$ & $4.11\times10^{-3}$ & $3.59\times10^{-3}$ & $4.24\times10^{-3}$ \\
        $E_{jump}$ & $4.75\times10^{-3}$ & $3.35\times10^{-3}$ & $2.11\times10^{-3}$ & $2.96\times10^{-3}$ \\
        $E_{flux}$ & $9.72\times10^{-3}$ & $4.49\times10^{-4}$ & $2.16\times10^{-4}$ & $3.63\times10^{-4}$ \\
		\noalign{\smallskip}\hline
	\end{tabular}
\end{table}

\begin{table}[H]
	\centering
	\caption{Relative $L^{2}$ errors for different numbers of modules and training points in Example \ref{E4}.}
	\label{E4_tab_module}
	\begin{tabular}{ccccc ccc}
		\hline\hline\noalign{\smallskip}	
    	Neurons & \ & 8 & 16 & 32  & 48\\
		\noalign{\smallskip}\hline\noalign{\smallskip}
        &$1$ & $7.57\times10^{-2}$ & $8.79\times10^{-3}$ & $4.11\times10^{-3}$ & $3.37\times10^{-3}$ \\
        Modules &$2$ & $4.34\times10^{-2}$ & $7.36\times10^{-3}$ & $4.77\times10^{-3}$ & $5.16\times10^{-3}$ \\
        &$3$ & $2.27\times10^{-2}$ & $7.19\times10^{-3}$ & $5.62\times10^{-3}$ & $4.84\times10^{-3}$ \\
		\noalign{\smallskip}\hline
	\end{tabular}
\end{table}

\begin{figure}[H]
    \centering
    \includegraphics[width=15cm,height=13cm]{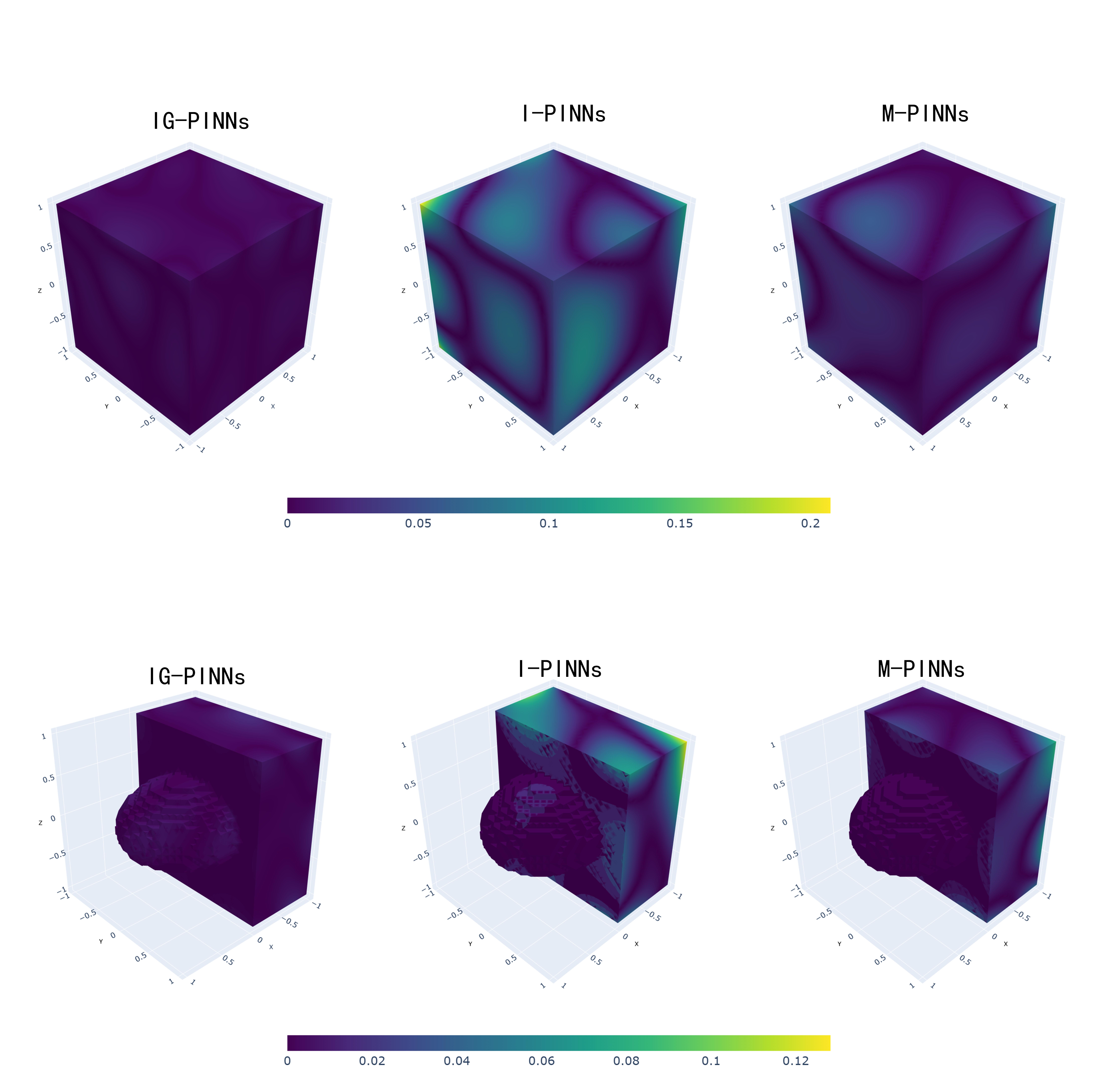}
    \caption{Example \ref{E4}, absolute errors for I-PINNs, M-PINNs and IG-PINNs. Top: the absolute error at the boundary of the computational domain. Bottom: the absolute error at the interface.}
    \label{E4_fig_error}
\end{figure}

\begin{figure}[H]
    \centering
    \includegraphics[width=14cm,height=6cm]{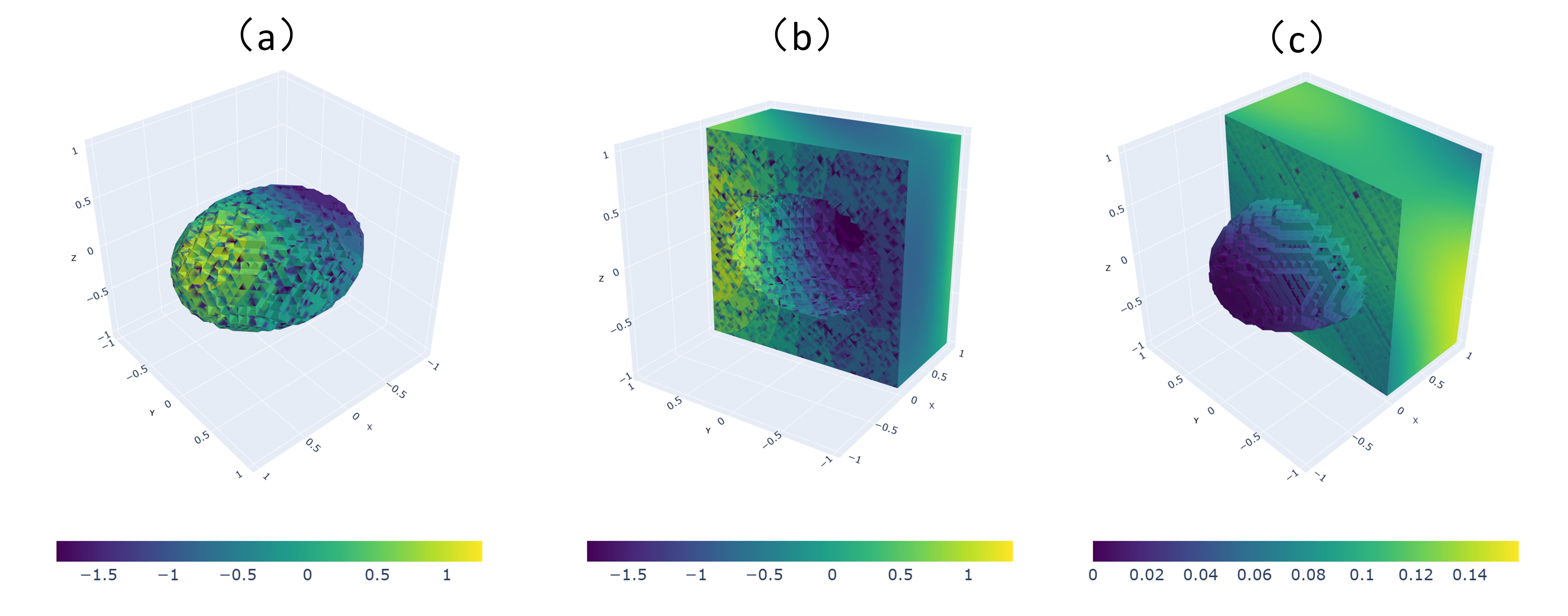}
    \caption{Example \ref{E4}, (\textbf{a}): Network component $u_{IG,1}(x,y,z)$. (\textbf{b}): Network component $u_{IG,2}(x,y,z)$. (\textbf{c}): Network component $\mu_{nn}(x,y,z)$. }
    \label{E4_fig_component}
\end{figure}

This test is similar to  Example \ref{E2}, \ref{E3}. The solution in subdomain $\Omega_{1}$ is only constrained by the interface conditions. IG-NNs contain some modules with interface information. Thus, IG-NNs explicitly focus on the interface, ensuring that the interface conditions are sufficiently optimized. Figure \ref{E4_fig_solution} shows the exact solution of problem (\ref{3.21}). Figure \ref{E4_fig_error} shows the absolute errors of IG-PINNs, I-PINNs, and M-PINNs at different locations, including interface and boundaries. The results show that the proposed method achieves the best approximation compared with I-PINNs and M-PINNs. Table \ref{E4_tab_error} shows the various errors of the IG-PINNs and the reference methods, including the maximum absolute error, the relative $L^2$ error. It can be seen that the errors of IG-PINNs are more than one magnitude lower than those of other algorithms. Figure \ref{E4_fig_component} shows the various network components of IG-PINNs. We test the performance of IG-PINNs at different numbers of modules and training points. Table \ref{E4_tab_module} provides a detailed comparison. The numerical results indicate that a smaller network scale can still achieve a good approximation, which avoids the use of larger network models. Table \ref{E4_tab_inter} provides a detailed error comparison for different numbers of interface points. A series of comparisons demonstrate the effectiveness of IG-PINNs.

\begin{example}\label{E5}
Finally, we consider a $3$D elliptic multiple interface problem in $\Omega=[-1,1]^{3}$.
The interfaces are as follows:

\begin{equation}\label{3.25}
\begin{split}
&\Gamma_{1}=\{(x,y,z):\phi_{1}(x,y,z)=0,\phi_{1}(x)=(x-0.25)^{2} + (y-0.25)^{2}+z^{2}-0.2^{2}\}, \\
&\Gamma_{2}=\{(x,y,z):\phi_{2}(x,y,z)=0,\phi_{2}(x)=(x-0.5)^{2} + (y+0.25)^{2}+z^{2}-0.15^{2}\},\\
&\Gamma_{3}=\{(x,y,z):\phi_{3}(x,y,z)=0,\phi_{3}(x)=x^{2} + (y+0.25)^{2}+(z-0.25)^{2}-0.1^{2}\}.\\
\end{split}
\end{equation}

The exact solution has the following form:
\begin{equation}\label{3.24}
u(x,y,z)=\left\{
	\begin{aligned}
	e^{xyz},    &\qquad(x,y,z)\in\Omega_{1},\\
	0.1sin(\pi(x+y))e^{(xyz)}, &\qquad(x,y,z)\in\Omega_{2},\\
    4x^{2} + y^{2} +2z^{2}, &\qquad(x,y,z)\in\Omega_{3},\\
    e^{x-z}cos(0.5\pi y),  &\qquad(x,y,z)\in\Omega_{4}.\\
	\end{aligned}
	\right.
\end{equation}
\end{example}

In this test, the computational domain is divided into four subdomains by the interfaces. We use a fully connected network with three hidden layers and four IG-NNs, each with one module, to approximate the exact solution. To highlight the performance of the proposed network framework, we choose I-PINNs and M-PINNs as reference algorithms. A total of 2,800 training points are used for training, including 1,000 interior points, 600 boundary points, and 400 interface points per interface. For the calculation of errors, 100 equally spaced nodes in each dimension (totaling $10^6$ test points) are used.

\begin{table}[H]
	\centering
	\caption{Comparison of errors in Example \ref{E5}.}
	\label{E5_tab_error}
	\begin{tabular}{ccccc ccc}
		\hline\hline\noalign{\smallskip}	
		 & I-PINNs & M-PINNs & IG-PINNs  \\
		\noalign{\smallskip}\hline\noalign{\smallskip}
		  $E_{M}$  & $1.05\times10^{-1}$ & $3.41\times10^{-2}$ & $1.01\times10^{-2}$ \\
          $E_{L^2}$ & $1.82\times10^{-2}$ & $5.37\times10^{-3}$ & $4.71\times10^{-4}$ \\
          $E_{jump,1}$ & $1.44\times10^{-2}$ & $7.14\times10^{-3}$ & $4.63\times10^{-3}$ \\
          $E_{flux,1}$ & $9.21\times10^{-2}$ & $6.88\times10^{-3}$ & $5.87\times10^{-3}$ \\
           $E_{jump,2}$ & $6.27\times10^{-2}$ & $1.18\times10^{-2}$ & $7.19\times10^{-3}$ \\
          $E_{flux,2}$ & $1.14\times10^{-2}$ & $4.33\times10^{-3}$ & $3.47\times10^{-4}$ \\
           $E_{jump,3}$ & $2.75\times10^{-2}$ & $8.17\times10^{-3}$ & $5.15\times10^{-3}$ \\
          $E_{flux,3}$ & $8.79\times10^{-3}$ & $5.62\times10^{-3}$ & $2.37\times10^{-3}$ \\
          Elapsed time & $387$(s) & $411$(s) & $724$(s) \\
		\noalign{\smallskip}\hline
	\end{tabular}
\end{table}

\begin{table}[H]
	\centering
	\caption{Comparison of relative $L^{2}$ errors in Example \ref{E5}.}
	\label{E5_tab_module}
	\begin{tabular}{ccccc ccc}
		\hline\hline\noalign{\smallskip}	
    	Neurons & \ & 8 & 16 & 32  & 48\\
		\noalign{\smallskip}\hline\noalign{\smallskip}
        &$1$ & $5.74\times10^{-3}$ & $4.71\times10^{-4}$ & $5.34\times10^{-4}$ & $5.79\times10^{-4}$ \\
        Modules &$2$ & $4.17\times10^{-3}$ & $5.77\times10^{-4}$ & $5.03\times10^{-4}$ & $6.35\times10^{-4}$ \\
        &$3$ & $2.28\times10^{-3}$ & $6.53\times10^{-4}$ & $4.97\times10^{-4}$ & $5.76\times10^{-4}$ \\
		\noalign{\smallskip}\hline
	\end{tabular}
\end{table}

\begin{figure}[H]
    \centering
    \includegraphics[width=15cm,height=8cm]{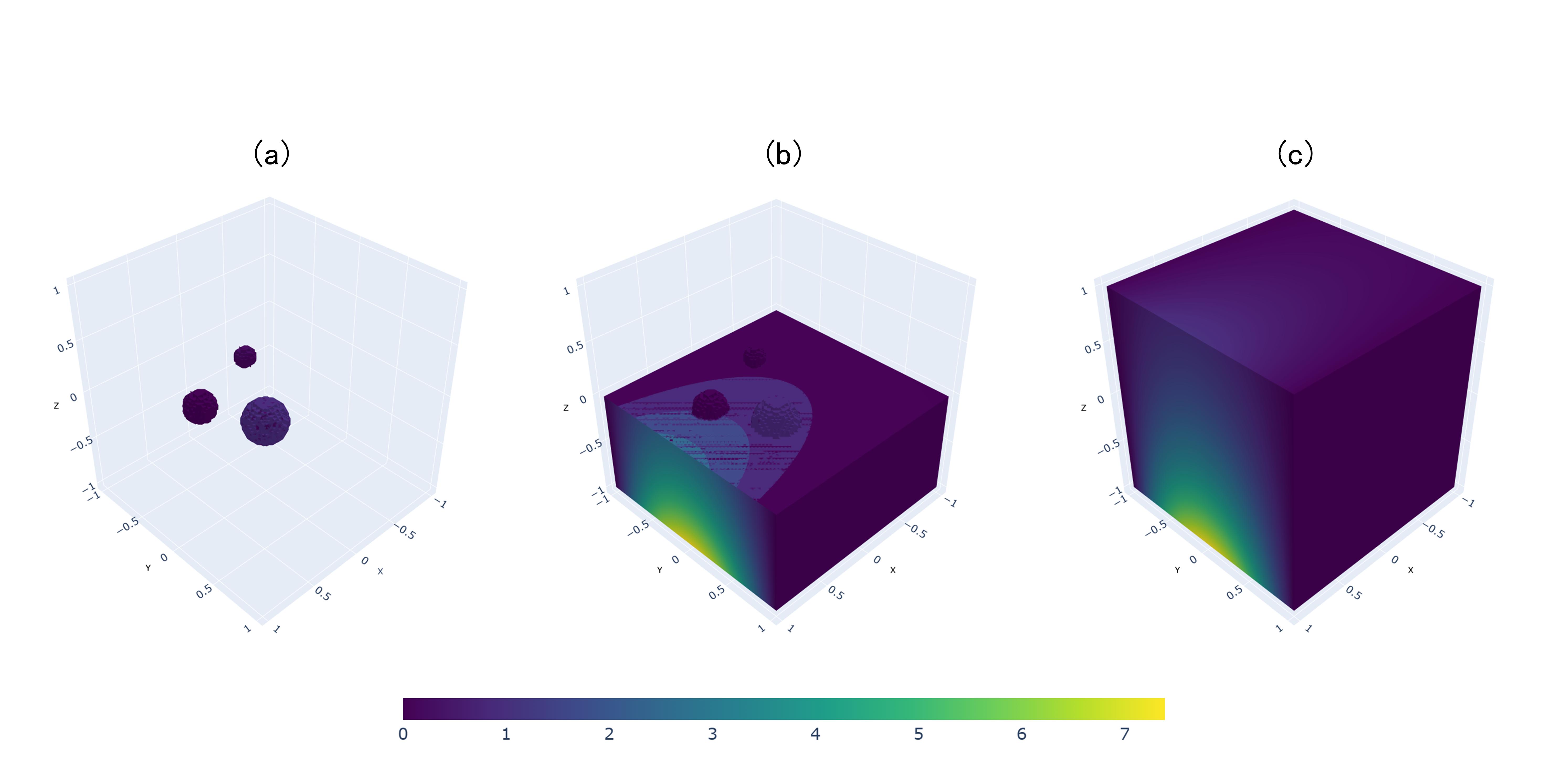}
    \caption{Example \ref{E5}, (\textbf{a}): The exact solution at the interface. (\textbf{b}): The exact solution inside the computational domain. (\textbf{c}): The exact solution at the boundary of the computational domain.}
    \label{E5_fig_solution}
    \end{figure}

\begin{figure}[H]
    \centering
    \includegraphics[width=15cm,height=19cm]{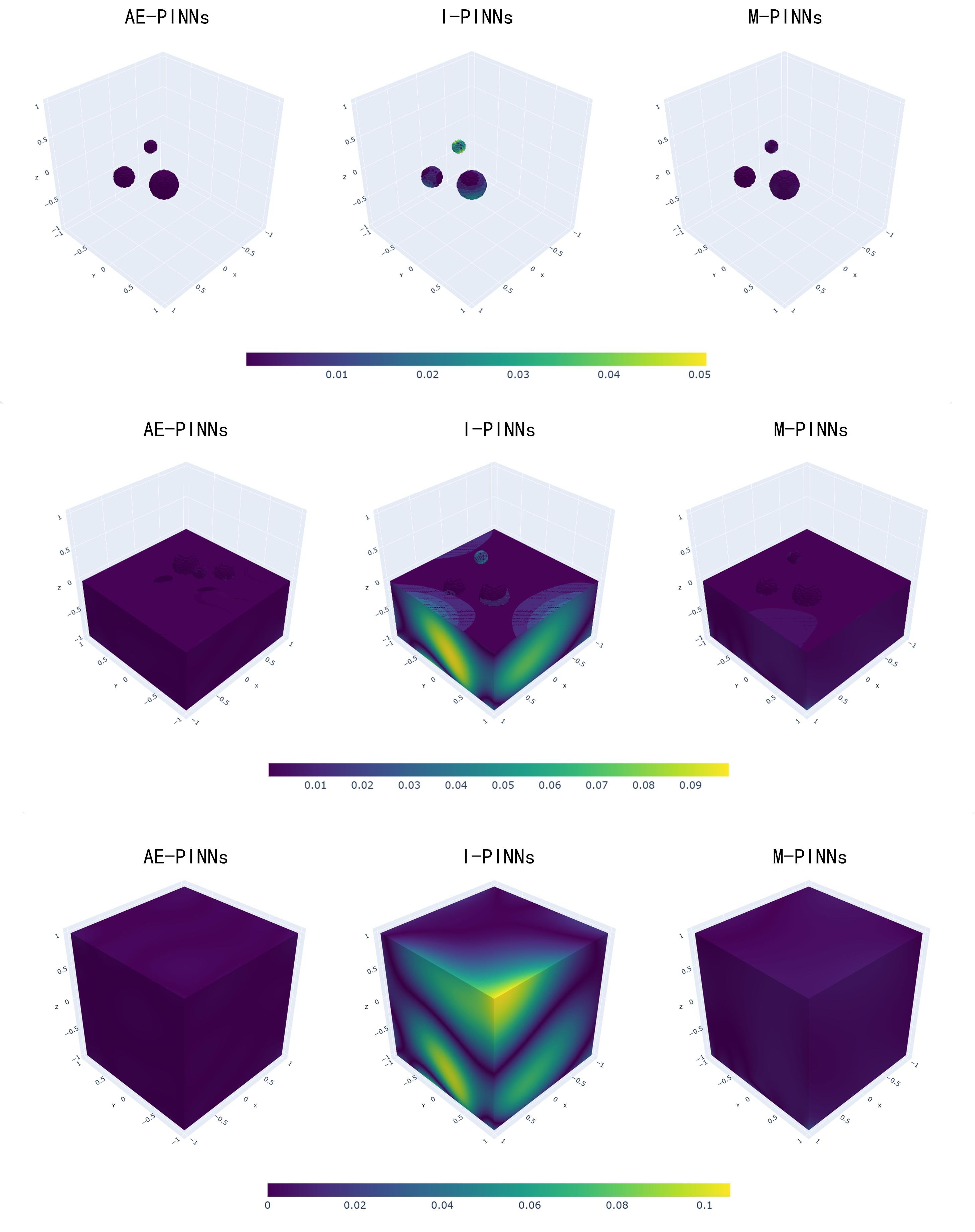}
    \caption{Example \ref{E5}, absolute errors for I-PINNs, M-PINNs and IG-PINNs. The
row-wise captions are; Top: the absolute error at the interface. Middle: the absolute error inside the computational domain. Bottom: the absolute error at the boundary of the computational domain.}
    \label{E5_fig_error}
    \end{figure}

\begin{figure}[H]
    \centering
    \includegraphics[width=14cm,height=6cm]{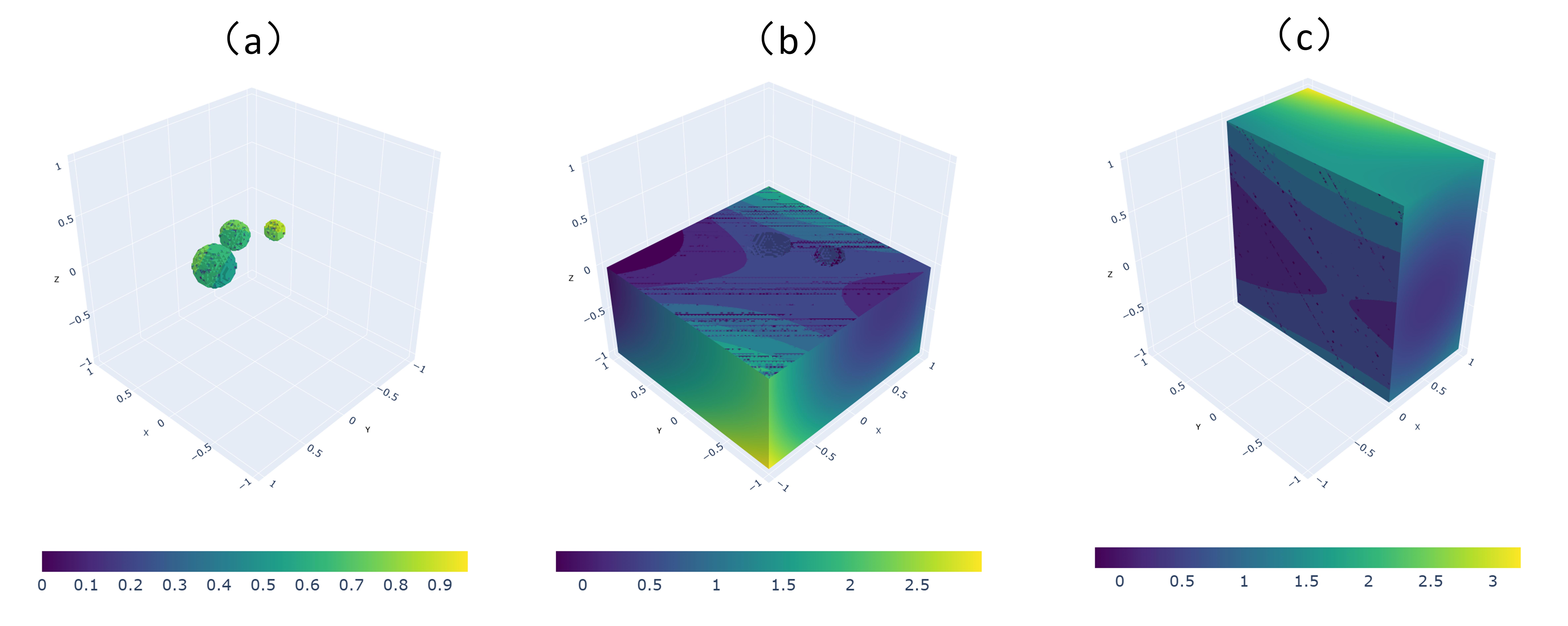}
    \caption{Example \ref{E5}, (\textbf{a}): Network component $u_{IG,i}(x,y,z), i=1,2,3$. (\textbf{b}): Network component $u_{IG,4}(x,y,z)$. (\textbf{c}): Network component $\mu_{nn}(x,y,z)$. }
    \label{E5_fig_com}
\end{figure}

Figure \ref{E5_fig_solution} shows the exact solution at different locations. Figure \ref{E5_fig_error} shows the absolute errors obtained by the three neural network algorithms, where the maximum absolute error is $1.01\times10^{-2}$ for IG-PINNs, $1.05\times10^{-1}$ and $3.41\times10^{-2}$ for I-PINNs and M-PINNs, respectively. Table \ref{E5_tab_error} provides a comparison of errors obtained by the three algorithms. Results show that the relative $L^{2}$ errors obtained by IG-PINNs reach $O(10^{-4})$, while those of I-PINNs and M-PINNs are in the range of $O(10^{-2})$, $O(10^{-3})$, respectively. Figure \ref{E5_fig_com} shows the various network components of IG-PINNs.

\begin{table}[htbp]
	\centering
	\caption{Comparison of errors in Example \ref{E5}.}
	\label{E5_tab_inter}
	\begin{tabular}{ccccc ccc}
		\hline\hline\noalign{\smallskip}	
    	Interface points & 20 & 50 & 100  & 200\\
		\noalign{\smallskip}\hline\noalign{\smallskip}
        $E_{L^2}$ & $5.13\times10^{-3}$ & $4.71\times10^{-4}$ & $5.63\times10^{-4}$ & $4.63\times10^{-4}$ \\
        $E_{jump,1}$ & $4.75\times10^{-3}$ & $4.63\times10^{-3}$ & $5.03\times10^{-3}$ & $4.29\times10^{-3}$ \\
        $E_{flux,1}$ & $9.72\times10^{-3}$ & $5.81\times10^{-3}$ & $6.19\times10^{-3}$ & $5.44\times10^{-3}$ \\
        $E_{jump,2}$ & $4.75\times10^{-3}$ & $7.19\times10^{-3}$ & $5.77\times10^{-3}$ & $6.21\times10^{-3}$ \\
        $E_{flux,2}$ & $9.72\times10^{-3}$ & $3.47\times10^{-4}$ & $4.41\times10^{-4}$ & $5.71\times10^{-4}$ \\
        $E_{jump,3}$ & $4.75\times10^{-3}$ & $5.15\times10^{-3}$ & $7.31\times10^{-3}$ & $5.83\times10^{-3}$ \\
        $E_{flux,3}$ & $9.72\times10^{-3}$ & $2.37\times10^{-4}$ & $2.73\times10^{-3}$ & $3.47\times10^{-3}$ \\
		\noalign{\smallskip}\hline
	\end{tabular}
\end{table}
 
We test the performance of IG-PINNs at different numbers of modules and training points. Table \ref{E5_tab_module} provides a detailed comparison. Numerical results indicate that IG-PINNs with a module can achieve a good approximation. Furthermore, we tested the model performance with different numbers of interface points, and the numerical results are presented in Table \ref{E5_tab_inter}. Combined with the previous examples, it can be seen that the IG-PINNs are sufficiently capable of solving interface problems compared to the I-PINNs and M-PINNs. However, it should be noted that IG-PINNs may require more training time compared to I-PINNs and M-PINNs.

\section{Summary and discussion}\label{section-4}
PINNs with continuous activation functions fail to accurately capture the discontinuities of solutions across the interfaces. To alleviate this limitation, we propose interface-gated physics-informed neural networks (IG-PINNs) for solving elliptic interface equations. The proposed algorithm involves the use of domain decomposition and coupling neural networks (we use two neural networks to approximate solution in each subdomain). Concretely, we use a fully connected neural network to capture the smooth behavior across the entire domain. In each subdomain separated by the interface, an interface-gated network is utilized to provide corrections at the interface.

The key contribution of this work is the proposal of interface-gated neural networks (IG-NNs). IG-NNs explicitly focus on the interface, ensuring that the interface conditions are optimized sufficiently during the training process. The performance of four neural network algorithms (PINNs, I-PINNs, M-PINNs and IG-PINNs) for solving interface problems is explored in detail through numerical examples. The results show that the PINNs provide a very bad approximation, which demonstrates the limitations of the PINNs for solving interface problems. Note that the relative $L^{2}$ error of the IG-PINNs is lower than that of the I-PINNs and M-PINNs by at least an order of magnitude for all examples, which demonstrates the effectiveness of the proposed method.

In this work, IG-PINNs result in higher computational cost compared with I-PINNs and M-PINNs. As a future extension of our research, we will focus on developing a simple and efficient network model. In addition, how to embed interface information into the network structure needs to be further discussed.\\

\section*{Declaration of competing interest}
The authors declare that they have no known competing financial interests or personal relationships that could have appeared to influence the work reported in this paper.

\section*{Acknowledgments}
The authors appreciate the constructive suggestions provided by the anonymous reviewers, which have significantly improved the quality of this paper.

This work was supported by NSFC Project (12431014) and Project of Scientific Research Fund of the Hunan Provincial Science and Technology Department (2024ZL5017).

\end{document}